\documentclass{rob}%

\usepackage{amssymb}
\usepackage{amsmath}
\usepackage[super]{cite}

\usepackage{array}
\usepackage{subfig}

\usepackage{multirow}
\usepackage{xcolor}
\usepackage{inputenc}
\usepackage{tikz}
\usetikzlibrary{shapes,arrows,calc,positioning}
\usepackage[percent]{overpic}
\usepackage{makecell}
\usepackage{boldline}
\usepackage{booktabs}
\usepackage[flushleft]{threeparttable}
\usepackage{color,soul}

\doi{10.1017/xxxx}
\volume{xx}
\issue{x}
\pubyear{2018}
\cpyear{2018}
\endpage{18}

\begin{document}
	
	\title[Time-Optimal Trajectory of Redundantly Actuated Robots]{A Computationally Efficient Algorithm to Find Time-Optimal Trajectory of Redundantly Actuated Robots Moving on a Specified Path}
	
	\author{Saeed Mansouri$\dagger$ \thanks{Corresponding author. E-mail: s\_mansouri@mech.sharif.edu}, Mohammad Jafar Sadigh$\ddagger$ and Masoud Fazeli\S}
	\affil{$\dagger$School of Mechanical Engineering, Sharif University of Technology, Tehran, Iran\\
		$\ddagger$School of Mechanical Engineering, College of Engineering, University of Tehran, Tehran, Iran\\
		{\S}Department of Mechanical Engineering, Isfahan University of Technology, Isfahan, Iran}
	
	\ADaccepted{MONTH DAY, YEAR. First published online: MONTH DAY, YEAR}
	
	\maketitle
	
	\begin{summary}
		\hl{A time-optimal} problem for redundantly actuated robots moving on a specified path is a challenging problem. Although the problem is well explored and there are proposed solutions based on phase plane analysis, there are still several \hl{unresolved} issues regarding calculation of solution \hl{curves}. In this paper, we explore the characteristics of \hl{the} maximum velocity curve (MVC) and propose an efficient algorithm to establish the solution curve. Then we propose a straightforward method to calculate the \hl{maximum or minimum} possible acceleration on the path based on the pattern of saturated actuators, which substantially reduces the computational cost. Two numerical examples are provided to illustrate the issues and the solutions.
	\end{summary}
	
	
\section{Introduction}
\label{sec:Cmpss}
Time-optimal motion planning remains an ongoing challenge for redundantly actuated robots, such as cooperative multi-manipulator systems, \hl{or} humanoid robots during double support phase. While several approaches have been proposed to \hl{tackle this issue}, large computational cost \hl{has yet remained} as a major obstacle. \hl{Here, we} propose a computationally efficient algorithm to address the time-optimal \hl{planning} problem for \hl{robots with actuation redundancy}. 
	
Several \hl{studies} have \hl{examined} the problem of minimum time motion \hl{for serial and parallel} manipulators \hl{while operating} on a specified path.
In general, there are two approaches \hl{to deal with such problem, including, the} Pontryagin maximum principle, and convex optimization. 
\hl{In the former approach, a system that moves in a specified path is regarded as having} only one degree of freedom (DOF). \hl{Accordingly, it is proven that a bang-bang solution,} in terms of the acceleration of the tip along the path, $ \ddot s $, \hl{accounts for the problem.}
\hl{This means that} the solution curve in \hl{the} $ s-\dot s $ plane could be obtained by \hl{successively integrating} the maximum and minimum values of $ \ddot s $ \cite{Bobrow1985}. \hl{This method, mainly known as the numerical integration method} \cite{Pham2014}, \hl{has been widely used in previous studies for time-optimal solutions in both normally actuated} \cite{Bobrow1985, Pfeiffer1987, Zlajpah1996, Kunz2012, Nguyen2016 , Shen2017 , Behzadipour2006} \hl{and over-actuated} \cite{Moon1990, Bobrow1990,  Moon1991, Moon1997, McCarthy1992, Ghasemi2008, Sadigh2013, Pham2012, Caron2017} \hl{robots. Guaranteed convergence of solution, low computational cost, and providing the solution directly in the trajectory space are of main advantages to this approach} \cite{Pham2015}.
	
In \hl{the} convex optimization approach, the time-optimal path tracking problem is transformed into a convex optimal control problem by discretization \hl{the path} $ s $ into $ N $ segments \cite{Verscheure2009, Zhang2016, Zhao2017}. Recently, this method has been used to find the minimum time motion of humanoid robots in multi-contact tasks \cite{Hauser2014}.
\hl{The} convex optimization method is \hl{indeed a more generic approach, as it not only applies to finding time-optimal solutions, but can also be used} to minimize any other index function.
\hl{Moreover, there are available commercial packages that exclusively deal with the convex optimization problem, making it a more accessible approach.}
	
\hl{In this paper, we mainly focus on the direct integration method for the time-optimal solution of redundantly actuated robots, while moving on a specified path. We aim to try and further reduce the computational costs associated with this method.} One should note that parallel manipulators are an important example of such systems for which the proposed method can be used.
	
In \hl{the} numerical integration method, the solution \hl{process} consists of two \hl{parts:} finding \hl{the} extremum acceleration in each step of solution, and extracting \hl{the} switching points.
A primary approach to \hl{finding} the extremum acceleration was first proposed in [\citen{Bobrow1985}] for non-redundant serial manipulators. The admissible range of $ \ddot s $ is determined at each ($ s $, $ \dot s $) under actuation constraints.
A more rigorous approach to find the maximum or minimum $ \ddot s $ was proposed in [\citen{Pfeiffer1987}] by taking \hl{advantage} of geometrical concepts.
In the latter method, inequality constraints on actuators are converted into a convex polygon in $ {\dot s}^2-\ddot s $ plane.
The solution provided in this study, however, was applied in serial robots with no redundant actuators.
Recently, this concept has been extended to over-actuated robots based on polytope projection technique \cite{Pham2015}.
\hl{Alternatively, the linear} programming (LP) method \cite{Moon1990, Bobrow1990, Moon1991, Moon1997} \hl{has also been used} to determine extremum acceleration for redundantly actuated systems. LP, \hl{however,} is a search based algorithm, which makes it numerically inefficient\hl{, and leads to large computational costs in the solution process.}
	
\hl{The first issue addressed in this paper is how to reduce the computational cost of calculating maximum/minimum acceleration along the path. To this end, we propose a method based on the pattern of saturated actuators. 
According to this method, }the saturation pattern of actuators is preserved in intervals of $ s $. \hl{It, hence,} reduces the challenges of solving LP \cite{Moon1990, Bobrow1990, Moon1991, Moon1997} or forming polygon \cite{Pham2015} at each step of integration, to a much simpler problem of finding the pattern of saturated actuators at some specific $ s $. \hl{It then follows by solving} a linear set of equations of motion at each step.
	
The second significant problem is to find \hl{the} switching points in cases \hl{where} more than one \hl{switching point} exists. Such \hl{problem occurs} when the velocity of the manipulator along the \hl{path exceeds} a maximum feasible value. In such cases, one needs to establish the maximum velocity curve (MVC) and find possible switching points \hl{on this curve,} known as critical points \cite{Pfeiffer1987}. 
One approach to \hl{find} MVC is \hl{to search} for the maximum feasible $ \dot s $ at each point of the path. To this end, one needs to find $ \dot s $ for which $ \ddot s_{max} = \ddot s_{min} $. This \hl{technique has been previously} applied for both serial and parallel manipulators in [\citen{Bobrow1985}] and [\citen{Moon1990, Bobrow1990, Moon1991, Moon1997}], respectively. \hl{However, it is evident that this method is numerically inefficient due to successive calculation of $ \ddot s_{max} $ and $ \ddot s_{min} $.} 
To reduce the computational cost, a method for direct calculation of MVC was suggested in [\citen{Pfeiffer1987}] for non-redundant serial manipulators.
Recently, the same technique has \hl{also} been used to determine MVC for over-actuated robots \cite{Pham2015}. Nevertheless, the construction of MVC for whole domain of $ s $ is still costly and numerically expensive.
	
Among the previous studies, three different strategies \hl{have been} proposed to determine the critical points. The first strategy that theoretically made an important advancement is based on the shooting method \cite{Bobrow1985}. Here, critical points are determined by looking for a solution curve which comes in contact with the MVC without crossing it. However, from \hl{a} practical point of view this approach was very difficult to apply and numerically inefficient. Another strategy which \hl{was} introduced in [\citen{Pfeiffer1987}] determines critical points by comparing the slope of MVC with the slope of solution curve for all \hl{points} located on MVC.
\hl{In both of the above-mentioned strategies, the construction of the MVC is necessary, which adds to the computational demands of the problem. This is while the MVC curve is not, per se, a solution curve, and nor does it play an important role in finding the switching points. What is crucial for finding the solution, however, is the lower boundary of the trapped area, as well as the locked area} \cite{Zlajpah1996}, \hl{which are both located under the MVC. This lower boundary, which itself is a solution trajectory, is in fact the switching curve. It seems, therefore, that it would be more efficient if one tries to construct the switching curve directly while skipping the construction of MVC as long as it is possible.}
	
\hl{In this direction, a straightforward approach has been proposed in} [\citen{Ghasemi2008}] \hl{to find  jump points on the MVC (a subset of critical points that are also known as zero-inertia points} \cite{Shiller1992, Pfeiffer1987, Pham2014}), \hl{to establish the switching curve based on that.} Such advancement could substantially reduce the computational cost. \hl{In some cases, however, such as finding the time-optimal solution for a} humanoid robot during double support phase \cite{Sadigh2013}, there remains a part of MVC whose construction is inevitable. \hl{Such a case motivated this study to further examine this issue and try to enhance the method proposed in} [\citen{Ghasemi2008}]. 
In this method, we introduce a new strategy to find the critical points by constructing MVC for minimum domain of $ s $. This strategy is based on direct construction of solution curve, taking advantage of analytical calculation of zero-inertia points\hl{, and to} find the domain of $ s $, if any, for witch no value of $ \dot s $ on solution curve is found during previous step. Construction of MVC is only necessary for such value of $ s $. 
We developed a new method to convert actuator constraints into a convex polygon in $ {\dot s}^2-\ddot s $ plane for over-actuated robots. 
The proposed method which directly constructs the polygon reduces the computational cost compared with the method proposed in [\citen{Pham2015}], which is based on iteration.
	
To summarize, in this study we tried to improve computational efficiency of solutions based on integration method for minimum time trajectory planning. \hl{The proposed approach is used} for cooperative multi-manipulator system (CMMS) with redundant actuators, moving on a specified path. \hl{Three} main contributions are made which are described in section \ref{sec:OM}. 
In section \ref{subsec:SPA}, an efficient method is proposed for computation of $ \ddot{s}_{max} $ and $ \ddot{s}_{min} $, which is to be found in each step of integration. 
In section  \ref{subsec:CP}, \hl{we discuss} which part of MVC \hl{could be skipped} by direct construction of solution curve, and which part \hl{is necessary to be constructed. The final} contribution \hl{in} this study is \hl{explained} in section \ref{subsec:POC}, \hl{where we describe} a novel method to directly establish the convex polygon of constraints for redundantly actuated robot \hl{as an} essential tool to construct MVC. This also helps to find pattern of saturated actuators.
\hl{Our proposed method is applied to solving two examples presented in section} \ref{sec:Exm}, \hl{and the elapsed computation time is compared with one of presented methods in} [\citen{Pham2015}]. The results of this comparison are given in section \ref{sec:EE}.

\section{Time Optimal Problem}

\begin{figure}[!t]
	\centering
	\includegraphics[width=0.9\linewidth]{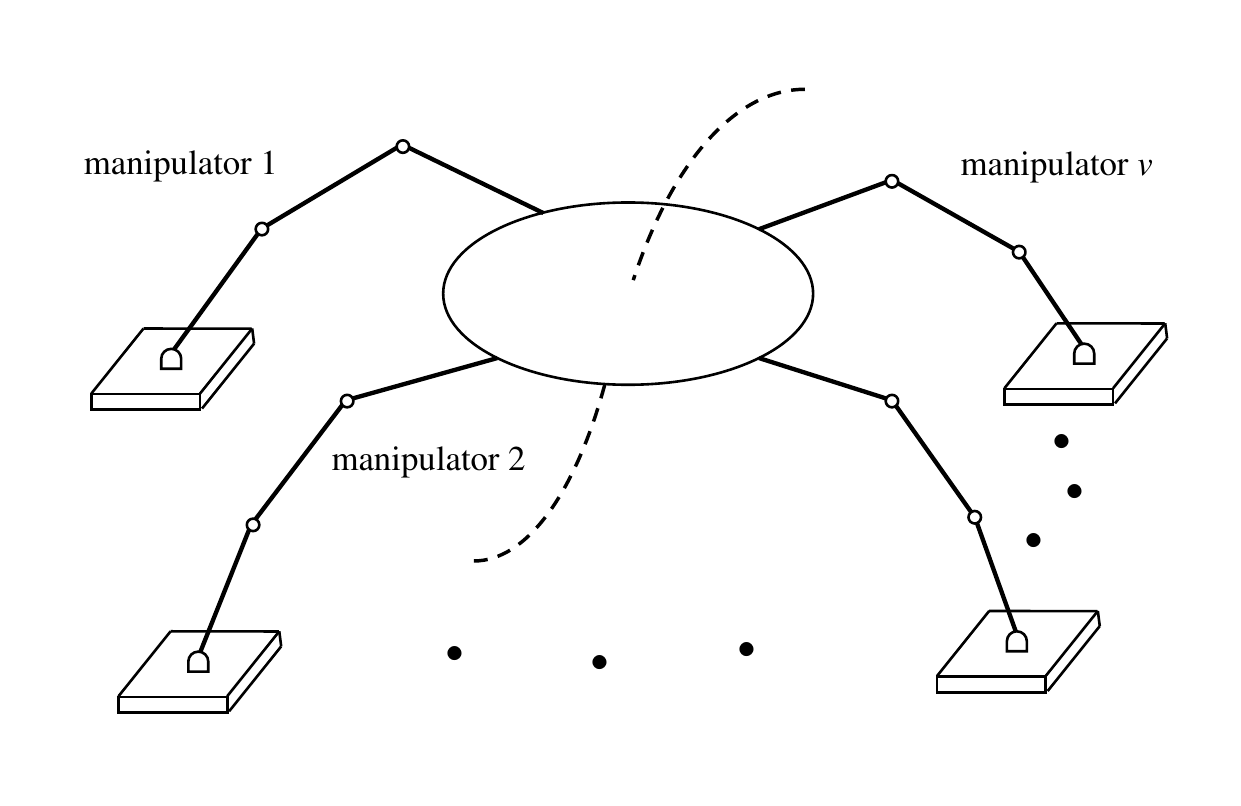}
	\caption{Schematic diagram for cooperative multi-manipulator system.}
	\label{fig:manipulator}
\end{figure}

An important set of redundantly actuated robots is CMMS which is considered to solved as an example in this paper.
Assume that the CMMS shown in Fig. \ref{fig:manipulator} is supposed to move a payload from an initial point to a final point on a specified path in minimum time subject to the actuator's limits. To state the time-optimal problem let us consider a CMMS composed of $ v $ similar, rigid non-redundant serial manipulators working on a single object. The motion of the payload in task space is defined by $ r $ coordinates; $ \boldsymbol x=[x_1, \dotso ,x_r]^T \in \mathbb{R}^{r}$ and the motion of the system is defined with $ n $ variables, $ q_1, \dotso ,q_n $ where $ \boldsymbol q=[\theta_{11}, \dotso ,\theta_{1r}, \dotso ,\theta_{v1}, \dotso ,\theta_{vr}]^T \in \mathbb{R}^{vr} $ and $ \theta_{ij} $  denotes the $ j $th joint coordinates of the $ i $th manipulator. Due to the kinematic chain configuration, the system is subject to $ p $ constraints. 
The reduced form of dynamic equations of such system can be written as

\begin{equation}
\label{eq:rde}
\mathbf{M}(\boldsymbol q)\ddot{\boldsymbol q}+\boldsymbol{h}(\boldsymbol q,\dot{\boldsymbol q})+\boldsymbol{g}(\boldsymbol q)=\mathbf{B}(\boldsymbol q) \boldsymbol{\tau}
\end{equation}

\noindent
where $ \mathbf{M}(\boldsymbol q)\in \mathbb{R}^{n-p\times n}$ is inertia matrix, $ \boldsymbol{h}(\boldsymbol q,\dot{\boldsymbol q})\in \mathbb{R}^{n-p}$ represents coriolis and centrifugal terms,  $ \boldsymbol{g}(\boldsymbol q)\in \mathbb{R}^{n-p}$ indicates gravitational effects, $ \mathbf{B}(\boldsymbol q)\in \mathbb{R}^{n-p\times m}$ is actuators coefficient matrix and $ \boldsymbol{\tau}\in \mathbb{R}^{m} $ is input torques \cite{Ghasemi2008}.

A desired path can be defined by

\begin{equation}
\label{eq:f}
x_{i}=f_{i}(s)\qquad
i=1, \dotso , r
\end{equation}

\noindent
in which $ x_i $ is one of the task space coordinates of manipulators and $ s $ is an independent non-dimensional parameter that uniquely describes the position of the system along the path. Also, the value of $ r $ is to be equal to $ n-p $ due to the fact that each manipulator is kinematically non-redundant.

We may state the time-optimal problem as:\medskip

$ Problem \ I: $ Find the optimum path $ \boldsymbol{q}^*(t) $ to minimize the elapsed time $ \Gamma=\int_{t_{i}}^{t_{f}} dt $ subject to

\begin{equation}
\label{eq:P1}
\begin{aligned}
&{\mathbf{M}}(\boldsymbol q)\ddot{\boldsymbol q}+{\boldsymbol h}(\boldsymbol q,\dot{\boldsymbol q})+{\boldsymbol g}(\boldsymbol q)={\mathbf{B}}(\boldsymbol q) \boldsymbol{\tau}\\
&\boldsymbol x=\boldsymbol f(s)\\
&\boldsymbol \tau_{min}\leq \boldsymbol \tau \leq  \boldsymbol \tau_{max}
\end{aligned}
\end{equation}

\noindent
with initial and final condition of

\begin{equation}
\begin{aligned}
&\boldsymbol q(t_i)=\boldsymbol  q_i,\qquad \dot{\boldsymbol{q}}(t_i)=\dot{\boldsymbol{q}}_i\\
&\boldsymbol q(t_f)=\boldsymbol  q_f,\qquad \dot{\boldsymbol{q}}(t_f)=\dot{\boldsymbol{q}}_f.\nonumber
\end{aligned}
\end{equation}

To solve $ problem \ I $ we restate it in terms of non-dimensional path parameter $ s $. To this end, one must substitute for $ \boldsymbol q $, $ \dot{\boldsymbol q} $ and $ \ddot{\boldsymbol q} $ in terms of $ s $, $ \dot{s} $ and $ \ddot{s} $ in (\ref{eq:rde}) to get

\begin{equation}
\label{eq:bar}
{\boldsymbol c}(s)\ddot s+{\boldsymbol d}(s)\dot s^2+{\boldsymbol e}(s)={\mathbf B}(s)\boldsymbol \tau
\end{equation}

\noindent
in which all $ {\boldsymbol c} $, $ {\boldsymbol d} $ and $ {\boldsymbol e} $ are members of $ \mathbb{R}^{n-p}$ \cite{Ghasemi2008}. 

The relation (\ref{eq:bar}) represents $ n-p $ equations with two states $ [s, \dot s] $. Any motion of the system which follows the prescribed path must satisfy all above $ n-p $ equations. Now, the optimization problem can be rewritten as:\medskip

$ Problem \ II: $ Find the optimum path $ s^*(t) $ to minimize the elapsed time $ \Gamma=\int_{t_{i}}^{t_{f}} dt $ subject to

\begin{equation}
\label{eq:P2}
\begin{aligned}
&{\boldsymbol c}(s)\ddot s+{\boldsymbol d}(s)\dot s^2+{\boldsymbol e}(s)={\mathbf B}(s)\boldsymbol \tau \\
&\boldsymbol \tau_{min}\leq \boldsymbol \tau \leq  \boldsymbol \tau_{max}
\end{aligned}
\end{equation}

\noindent
with initial and final condition of

\begin{equation}
\begin{aligned}
&s(t_i)=0,\qquad {\dot s}(t_i)=\dot s_i \\
&s(t_f)=1,\qquad {\dot s}(t_f)=\dot s_f. \nonumber
\end{aligned}
\end{equation}

As we stated earlier in introduction, finding the solution of $ problem \ II $ is equivalent to solve two major problems of finding extremum acceleration and extracting switching points.
In the next section, after a detailed review of some background, we propose an efficient method to calculate the extremum acceleration based on pattern of saturated actuators.

\section{Optimization Method}
\label{sec:OM}

\subsection{Saturation Pattern of Actuators}
\label{subsec:SPA}

The form of  $ problem \ II $ resembles that of the time-optimal for non-redundant serial manipulators with \hl{the} main difference \hl{being in} matrix $ \mathbf B $. In case of systems with no redundant actuators, this matrix is invertible and could be reduced into a unity matrix; which means that 
only one actuator \hl{does appear} in each equation.  
However, for CMMS which has redundancy in actuators, $ \mathbf{{B}} $ is a non-square matrix, which means that there \hl{are always} more than one actuator torques in each equation of $ problem \ II $. That is why \hl{the} provided methods  for serial manipulators in [\citen{Bobrow1985}] and [\citen{Pfeiffer1987}] cannot be used to calculate $  {\ddot s}_{max} $ or $  {\ddot s}_{min} $ for CMMS.
This issue was independently addressed in [\citen{Moon1990}] and [\citen{Bobrow1990}], where the problem of finding extremum acceleration was introduced as:

\medskip
$ Problem \ III: $ Find $ \boldsymbol \tau $ to maximize or minimize $ \ddot s $ subject to

\begin{equation}
\label{eq:P3}
\begin{aligned}
&{\boldsymbol c}(s)\ddot s+{\boldsymbol d}(s)\dot s^2+{\boldsymbol e}(s)={\mathbf B}(s)\boldsymbol \tau \\
&\boldsymbol \tau_{min}\leq \boldsymbol \tau \leq  \boldsymbol \tau_{max}.
\end{aligned}
\end{equation}

This problem could be solved using a linear programming (LP) strategy. One should note that this problem must be solved for all values of $ s $ and $ \dot s $ on \hl{the} solution curve, which, clearly, \hl{requires a huge amount of computation.}

It is demonstrated in [\citen{McCarthy1992}] that for a robotic system with $ n $ generalized coordinates, $ p $ holonomic constraints and $ m $ actuators, at least $ m-n+p+1 $ actuators are in \hl{the} saturated state. 
Therefore, solving $ problem \ III $ results in finding $ m-n+p+1 $ saturated actuators which \hl{leads} to $ \ddot s $ being maximized or minimized. The solution gives the value of $ \ddot s $ and $ \boldsymbol \tau $.  
Here, we introduce a method based on the saturation pattern of actuators, which reduces \hl{the} solution of $ problem \ III $ into \hl{the} solution of a set of linear equations of motion.

Considering the fact that when $ \ddot s $ accepts its extremum, at least $ m-n+p+1 $ actuators are saturated \cite{McCarthy1992}, one can easily compute non-saturated actuators and $ \ddot s $ from (\ref{eq:bar}). To this end, \hl{we will have}

\begin{equation}
\label{eq:E4}
{\boldsymbol c}(s)\ddot s+{\boldsymbol d}(s)\dot s^2+{\boldsymbol e}(s)=
\begin{bmatrix} {\mathbf B}_{ns} & {\mathbf B}_{s} \end{bmatrix}
\begin{bmatrix} \boldsymbol{\tau}_{ns} \\ \boldsymbol{\tau}_{s} \end{bmatrix}         
\end{equation}

\noindent
in which $ \boldsymbol \tau_{ns}\in \mathbb{R}^{n-p-1} $ and $ \boldsymbol \tau_{s}\in \mathbb{R}^{m-n+p+1} $.

For every point along the solution trajectory in $ s-\dot s $ plane, one can always \hl{find} a neighborhood in which \hl{the} saturation pattern of actuators does not change. Therefore, we may use a set of saturated actuators at one point to calculate $ \ddot s $ and $ \boldsymbol \tau_{ns} $ for other points in such \hl{a} neighborhood.
So the problem of determining extremum acceleration, $ problem \ III $, reduces to\medskip 

$ Problem \ IV: $ Assuming $ \boldsymbol {\tau}_s $ to be known, find $ \boldsymbol \tau_{ns} $ and $ \ddot s $ from

\begin{equation}
\label{eq:E7}
\begin{bmatrix} {\boldsymbol c}(s) & -{\mathbf B}_{ns} \end{bmatrix}
\begin{bmatrix} \ddot s \\ {\tau}_{ns} \end{bmatrix}  
=-{\boldsymbol d}(s)\dot s^2-{\boldsymbol e}(s)+{\mathbf B}_{s} \boldsymbol{\tau}_{s}.
\end{equation} 

Clearly\hl{, the} solution of (\ref{eq:E7}) which is a linear set is much easier and involves less computation compared \hl{to} (\ref{eq:P3}), which is a linear programming problem.
To solve the problem \hl{with} hand, first, the saturated actuators set which leads to \hl{maximizing or minimizing} $ \ddot s $ is determined for one point by constructing the polygon of constraints in $ {\dot s}^2-\ddot s $ plane as will be explained in section \ref{subsec:POC}, or by solving LP problem, i.e. $ problem \ III $, for that point. \hl{The next} point in $ s-\dot s $ plane is found by integration of equations of motion using this value of $ \ddot s $ and $ \boldsymbol{\tau} $. For \hl{the} next step, we do not need to solve the LP problem. \hl{Instead, because $ \boldsymbol{\tau}_{s} $ remains the same, we can} calculate $ \ddot s $ and $ \boldsymbol{\tau}_{ns} $ from (\ref{eq:E7}).

If none of the computed values for $ \boldsymbol{\tau}_{ns} $ \hl{exceeds the actuator's bound,} the solution is acceptable\hl{; and therefore,} we may proceed to \hl{the} next step. \hl{However,} if at least one of the values of $ \boldsymbol{\tau}_{ns} $ exceeds its bound, it means that the pattern of saturation is changed and we need to establish the polygon to find real set of saturated actuators.
Fig. \ref{fig:FD} shows the flow diagram of the proposed algorithm which we \hl{would call SPA hereafter.}

\tikzstyle{startstop} = [rectangle, rounded corners, minimum width=1.5cm, minimum height=0.5cm,text centered, draw=black, fill=red!30, rounded corners]
\tikzstyle{io} = [trapezium, trapezium left angle=70, trapezium right angle=110, minimum width=1.5cm, minimum height=0.5cm, text centered, draw=black, fill=blue!30, rounded corners]
\tikzstyle{process} = [rectangle, minimum width=1.5cm, minimum height=0.5cm, text centered, draw=black, fill=orange!30, rounded corners]
\tikzstyle{decision0} = [diamond, minimum width=0.5cm, minimum height=0.5cm, text centered, text width=1.4cm, draw=black, fill=green!30, rounded corners]
\tikzstyle{decision1} = [diamond, minimum width=1.5cm, minimum height=0.5cm, text centered, draw=black, fill=green!30, rounded corners]
\tikzstyle{decision2} = [diamond, minimum width=0.5cm, minimum height=0.5cm, text centered, text width=1.4cm, draw=black, fill=green!30, rounded corners]
\tikzstyle{arrow} = [thick,->,>=stealth]

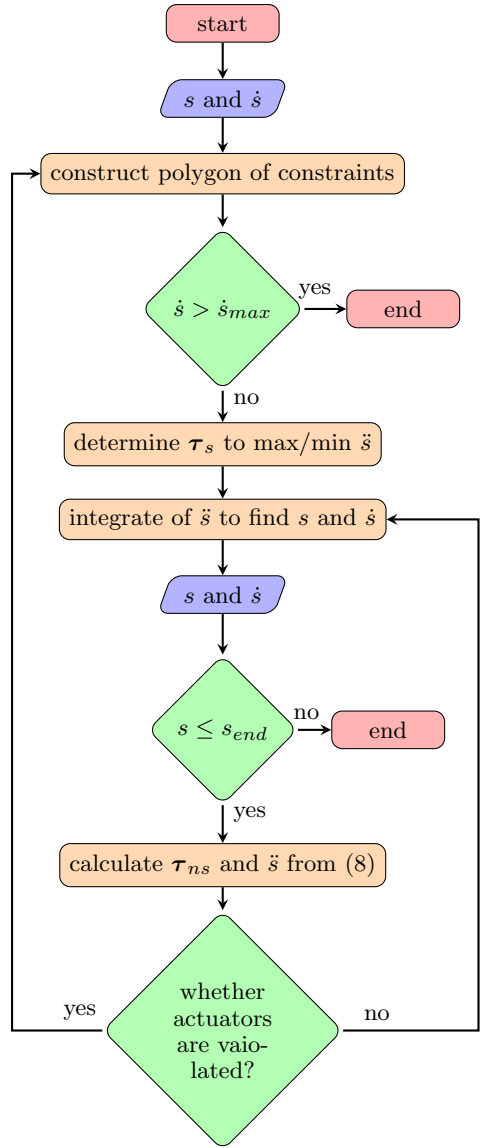
\begin{figure}[!t]
	\centering
	\begin{tikzpicture}[font=\footnotesize]
	
	\node (start) [startstop] {start};	
	\node (in1) [io, below of=start] {$ s $ and $ \dot s $};
	\node (pro2) [process, below of=in1] {construct polygon of constraints};
	\node (dec0) [decision0, below of=pro2, yshift=-0.8cm] {$ \dot{s}>\dot{s}_{max} $};
	\node (end0) [startstop, right of=dec0, xshift=1.4cm] {end};
	\node (pro4) [process, below of=dec0, yshift=-0.8cm] {determine $ \boldsymbol {\tau}_s $ to max/min $ \ddot s $};
	\node (pro5) [process, below of=pro4] {integrate of $ \ddot s $ to find $ s $ and $ \dot s $};
	\node (out1) [io, below of=pro5] {$ s $ and $ \dot{s} $};
	\node (dec3) [decision1, below of=out1, yshift=-0.8cm] {$ s \leq s_{end} $};
	\node (pro6) [process, below of=dec3, yshift=-0.8cm] {calculate $ \boldsymbol {\tau}_{ns} $ and $ \ddot s $ from (\ref{eq:E7})};
	\node (dec2) [decision2, below of=pro6, yshift=-1.2cm] {whether actuators are vaiolated?};		
	\node (end1) [startstop, right of=dec3, xshift=1.2cm] {end};
	
	\draw [arrow] (start) -- (in1);
	\draw [arrow] (in1) -- (pro2);
	\draw [arrow] (pro2) -- (dec0);
	\draw [arrow] (dec0) -- node[near start][above] {yes} (end0);
	\draw [arrow] (dec0) -- node[near start][right] {no} (pro4);	
	\draw [arrow] (pro4) -- (pro5);
	\draw [arrow] (pro5) -- (out1);
	\draw [arrow] (out1) -- (dec3);
	\draw [arrow] (dec3) -- node[near start][right]{yes}  (pro6);
	\draw [arrow] (dec3) -- node[near start][above]{no} (end1);
	\draw [arrow] (pro6) -- (dec2);
	\draw [arrow] (dec2) -- ($(dec2.west)+(-1.2,0)$) node[near start][above] {yes} |- (pro2);
	\draw [arrow] (dec2) -- ($(dec2.east)+(1.8,0)$) node[near start][above] {no} |- (pro5);
	
	\end{tikzpicture}
	\caption{Flow diagram of the SPA strategy}
	\label{fig:FD}
\end{figure}

By applying  this innovative method, while the $ problem \ IV $ is solved in each integration step, determining the pattern of saturated actuators which needs more computational effort, is done only in some limited points, where \hl{the} pattern of saturation changes.
This method, hence, significantly reduces the amount of computation.

\subsection{Critical Points}
\label{subsec:CP}
An important step to find the solution is determining the points to switch between \hl{the} minimum and maximum acceleration. For this purpose, the first stage is calculating the critical points which are \hl{candidates} for switching points on the MVC. 
The concept of critical points was first introduced in [\citen{Pfeiffer1987}].
Points on the MVC are categorized as sink or source depending on whether the slope of solution curve in phase plane, $ \ddot s/\dot s $, is greater or smaller than the slope of the MVC, $ k^\prime $.
Any solution starting from a source continues in feasible region, whereas any solution starting from a sink immediately leaves the feasible region. By critical points, we mean the points for which the nature of MVC is neither \hl{a sink  nor a source}.

In this section, we propose a new strategy to find the critical points. This method is an improvement of the method presented in [\citen{Ghasemi2008}]. 
This method involves two main parts. In the first part, we analytically determine zero-inertia points, i.e. the points which satisfy the following equation \cite{Ghasemi2008}

\begin{equation}
{{\overline c}}_i(s^\dagger)=0.
\label{eq:ci}
\end{equation}

\noindent
Then we construct the solution curve using these zero-inertia points. In the second part, we construct the portion of MVC \hl{has no points on the} solution curve with the same value of $ s $. This will be discussed in more details in section \ref{sec:SA} while describing the solution algorithm. One should note that 
not all $ s^\dagger $ which satisfy (\ref{eq:ci}) are feasible zero-inertia point. To check if $ s^\dagger $ is a feasible point, it is just enough to check whether $ \dot{s}^\dagger $ associated with $ s^\dagger $ violates actuator constraints or not.

One should also note that not all feasible zero-inertia points located on MVC are candidates for switching points. Depending on the change of the nature of boundary on two sides of a zero-inertia point, it might be categorized as one of \hl{the four} types of sink to source, source to sink, sink to sink and source to source. 
Out of these, only sink to source zero-inertia points are critical points. 
To realize the type of a zero-inertia point without constructing the MVC, one needs to establish the solution from such points in backward and forward directions with $ \ddot s_{min} $ and $ \ddot s_{max} $, respectively. It is done by integrating equations of motion backward and forward starting from that point. 
If the zero-inertia point is of the sink--sink type, for the forward direction, the \hl{actuator's bounds are immediately violated,} which means that \hl{the} forward path instantly enters the non-feasible region (NFR). 
Also, if the point is of the source--source type, the \hl{actuator's limits are immediately violated for the} backward direction. 
Similarly, if the point is of the source--sink type, for both the forward and backward \hl{directions actuator's bounds} are violated immediately. However, solutions starting from a sink--source point remain in \hl{the} feasible area.

In the developed strategy, in order to determine the switching points, only some parts of the MVC are constructed. It is more efficient than the other studies, e.g. [\citen{Pham2015}], which establish MVC in the whole domain of $ s $.

\subsection{Polygon of Constraints}
\label{subsec:POC}

As mentioned in \hl{the} previous section, there are points on the path for which we need to find $ \ddot s_{max} $ or $ \ddot s_{min} $ and the pattern of saturated actuators. In addition to that, one might also need to find the maximum possible $ \dot s $ at each $ s $ to establish MVC. One efficient method to solve both of the above mentioned problems is to establish \hl{the} polygon of constraints which was first reported in [\citen{Pfeiffer1987}] for serial manipulators.

To construct the polygon for CMMS systems, we attempt to extend the method proposed in [\citen{Pfeiffer1987}].
For this, consider relation (\ref{eq:bar}) which contains $ n-p $ equations. 
Choosing $ n-p $ independent actuators $ \boldsymbol \tau_{a}\in \mathbb{R}^{n-p} $, one might rewrite relation (\ref{eq:bar}) as

\begin{equation}
\label{eq:E4b}
{\boldsymbol c}(s)\ddot s+{\boldsymbol d}(s)\dot s^2+{\boldsymbol e}(s)=
\begin{bmatrix} {\mathbf B}_{a} & {\mathbf B}_{b} \end{bmatrix}
\begin{bmatrix} \boldsymbol{\tau}_{a} \\ \boldsymbol{\tau}_{b} \end{bmatrix}.   
\end{equation}

\noindent
Here by independent \hl{actuator,} we mean a set of them which are sufficient to statistically balance the system at \hl{each $ s $,} regardless of the choice of $ \boldsymbol \tau_{b} $. In above \hl{equation,} $ \boldsymbol \tau_{b}\in \mathbb{R}^{m-n+p} $ depicts the remaining actuators.
Premultiplying (\ref{eq:E4b}) by $ {{\mathbf B}_{a}}^{-1} $ leads to

\begin{equation}
\label{eq:E5}
{\overline{\boldsymbol c}}(s)\ddot s+{\overline{\boldsymbol d}}(s)\dot s^2+{\overline{\boldsymbol e}}(s)
=\boldsymbol{\tau}_{a} + {\overline{\mathbf  B}} \boldsymbol{\tau}_{b}
=\boldsymbol{\tau}^\prime         
\end{equation}

\noindent
in which $ {\overline{\boldsymbol c}}={{\mathbf B}_{a}}^{-1} {\boldsymbol c}\in \mathbb{R}^{n-p} $,
$ {\overline{\boldsymbol d}}={{\mathbf B}_{a}}^{-1} {\boldsymbol d}\in \mathbb{R}^{n-p} $,
$ {\overline{\boldsymbol e}}={{\mathbf B}_{a}}^{-1} {\boldsymbol e}\in \mathbb{R}^{n-p} $ and
$ {\overline{\mathbf  B}}={{\mathbf B}_{a}}^{-1} {{\mathbf B}}_b\in \mathbb{R}^{n-p\times m-n+p} $.
One should note that existence of $ {{\mathbf B}_{a}}^{-1} $ is guaranteed by the virtue of independency of actuators $ \boldsymbol{\tau}_{a} $.

In any equation of (\ref{eq:E5}), all $ {\boldsymbol \tau}_b $ and one of $ {{\tau}_{a}}_i $, which equal to $ m-(n-p)+1 $ \hl{actuators appear.} 
By selecting the appropriate values of $ \boldsymbol{\tau}_{b} $ and $ \tau_{{a}_i} $ from their lower or upper bound, one could obtain the minimum and maximum values for $ {\tau}^\prime_i $.
Therefore, each equation of (\ref{eq:E5}) shows two straight lines on the plane $ \dot s^2-\ddot s $. The area out of these two lines indicates a region in which for any $ (\dot s^2, \ddot s) $ at least one of the actuator's limits is violated so this point is non-feasible. 
The feasible area is restricted to points inside this polygon. Changing the choice of \hl{$ \boldsymbol \tau_b $,} other polygons are also generated.
The common area of all these polygons is considered as \hl{the} feasible area of CMMS. This area constitutes a unique polygon which works in the same way as the one established for serial manipulators \cite{Pfeiffer1987}. 

The polygon constructed for CMMS is the common area restricted by $ C(m, m-(n-p))\times n-p $ pairs of lines. Here, $ C $ depicts \hl{the} combination function. For any given $ (\dot s^2, \ddot s) $ inside the polygon, there is one feasible $ \boldsymbol \tau $ that satisfies the actuators saturation limits. On the edge of \hl{the} polygon, $ m-(n-p)+1 $ actuators are in saturation \hl{conditions.} $ m-(n-p) $ are those of $ \boldsymbol \tau_b $  and the other one is  the corresponding $ {{\tau}_{a}}_i $. 
At any vertex of \hl{ the polygon,} at least $ m-(n-p)+2 $ actuators are in saturation \hl{conditions.} The set of saturated actuators on each edge is the common ones at two vertices located at the ends of the edge. Therefore, by knowing the saturated actuators at vertices of \hl{the} polygon, the pattern of saturation on each edge is determined.

In this paper, we present a new method to establish the polygon for CMMS by finding its vertices.
By applying this method, in addition to the position of vertices in $ {\dot s}^2-\ddot s $ plane, the set of saturated actuators at each vertex is also obtained. Our proposed strategy to construct the polygon is as follows:

\begin{itemize}
	\item[$ i. $] Choose a state of $ \boldsymbol \tau_b $ and extract relation (\ref{eq:E5}) for it.
	
	\item[$ ii. $] For each equation of (\ref{eq:E5}), establish the corresponding pair of lines on $ \dot s^2-\ddot s $ plane.
	
	\item[$ iii. $] Calculate all intersection points of each two lines which \hl{have} the same value of $ \boldsymbol \tau_b $. 
	
	\item[$ iv. $] For each intersection point, compute the value of $ \boldsymbol \tau_a $  \hl{from} (\ref{eq:E5}) and check whether $ \boldsymbol \tau_a $ violates its saturation limits. If yes, this point could not be a vertex of polygon and should be ignored.
	
	\item[$ v. $] Continue step $ i $ to $ iv $ until all possible states of $ \boldsymbol \tau_b $ are selected. After this, all vertices of the polygon are obtained.
	
\end{itemize}

The developed method is a straightforward and efficient strategy to find the polygon.
It is different from the approach proposed in [\citen{Pham2015}] that is based on \hl{the} polytope projection technique. This technique iteratively attempts to find polygon by expanding the edges which leads to a trade-off between accuracy and computational cost.

\section{Solution Algorithm}
\label{sec:SA}

\begin{figure}[!t]
	\centering
	\includegraphics[trim=3cm 13.2cm 3cm 6cm, clip=true, width=1\linewidth]{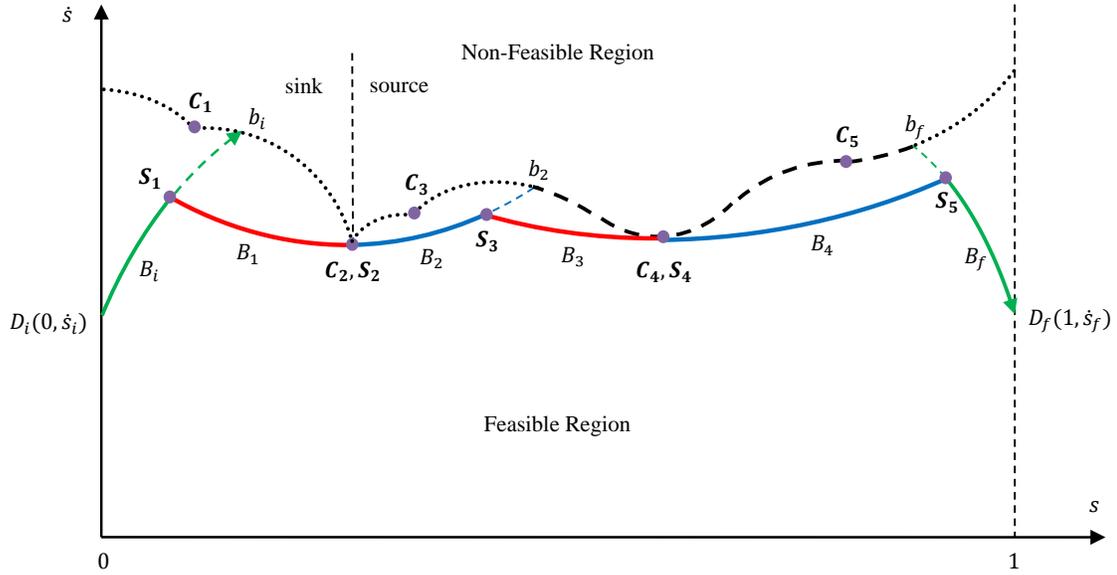}
	\caption[Schematic diagram of the construction of the solution curve.]{Schematic diagram of the construction of the solution curve. In this figure,
		$ D_i $ and $ D_f $ are initial and final points of solution curve,
		$ C_j $ is $ j $th critical point, 
		$ S_j $ is $ j $th switching point, 
		$ B_i $ and $ B_f $ are initial and final branches of solution curve, 
		$ B_j $ is $ j $th branch of solution curve, 
		$ b_i $ and $ b_f $ are points that initial and final solution curves enter NFR, and
		$ b_j $ is point that $ j $th solution curve enters NFR.
	}
	\label{fig:Schematic}
\end{figure}

In this section, we provide the final algorithm to solve the time-optimal problem of redundantly actuated robots moving on a specified path. The algorithm utilizes strategies introduced in previous sections to find critical points and determine the pattern of saturated actuators.
In phase plane, the higher solution trajectory has the shorter traveling time. By considering this fact, we propose the following algorithm to find the time-optimal trajectory:

$ Step \ 1: $
Using SPA strategy, start from initial point, $ D_i(0,\dot s_i) $, in forward direction with $ \ddot s_{max} $.
Similarly, start from final point, $ D_f(1,\dot s_f) $, in backward direction with $ \ddot s_{min} $. Continue these two procedures until either the line of $ s = 0 $ or $ s = 1 $ is crossed or \hl{actuator's bounds} are violated, i.e. solution trajectories enter the NFR (points $ b_i $ and $ b_f $ in Fig. \ref{fig:Schematic}).
If one or both of the curves \hl{intersect} the line of $ s = 0 $ or $ s = 1 $, solution trajectories cross each other at feasible region (FR). Therefore, the algorithm ends and the solution curve is obtained with one switching point. Otherwise, zero-inertia points should be determined, like the case shown in Fig. \ref{fig:Schematic}.

$ Step \ 2: $
Apply the strategy proposed in section \ref{subsec:CP} to calculate a set of zero-inertia points with the type of sink--source which are located on MVC between $ b_i $ and $ b_f $, e.g. $ C_2 $ and $ C_3 $ in Fig. \ref{fig:Schematic}. If there \hl{are} no such points skip to $ Step \ 6 $.

$ Step \ 3: $
Run SPA strategy by starting from the lowest point of the set, $ C_2 $, in backward and forward \hl{directions with $ \ddot s_{min} $ and $ \ddot s_{max} $,} respectively. This process is continued until either the solution trajectories generated in previous steps is crossed ($ S_1 $) or no acceptable solution fulfills the actuator's limits (point $ b_2 $ on MVC). 

$ Step \ 4: $
Ignore all other points of the set, $ C_3 $, which are located above the generated solution curve inside the inadmissible area, i.e. trapped and locked \hl{area introduced} in [\citen{Zlajpah1996}].

$ Step \ 5: $
Repeat $ Steps \ 3 $ and $ 4 $ from the lowest remaining points of the set. 
Continue the same process for other points of the set until either the generated solution curve covers the whole domain of $ s $ from $ 0 $ to $ 1 $ or there remains no zero-inertia point outside \hl{the} trapped and locked area, like the case shown in Fig. \ref{fig:Schematic}. If both of the initial and final solution trajectories which are established in $ Step \ 1 $ are crossed with the generated curve, the algorithm ends and the solution curve is obtained. Otherwise, some parts of MVC should be constructed, like the case shown in Fig. \ref{fig:Schematic}.

$ Step \ 6: $
Construct portions of MVC which are not covered by \hl{the} previously generated solution curve, i.e. dashed curve $ b_2b_f $. It is done by establishing polygons for such intervals using the algorithm proposed in \hl{Section} \ref{subsec:POC}.

$ Step \ 7: $
By calculating $ \ddot s/\dot s - k^\prime $, find a set of smooth critical points where the characteristics of the MVC constructed in $ Step \ 6 $ \hl{changes} from sink to source, e.g. points $ C_4 $ and $ C_5 $.

$ Step \ 8: $
Repeat $ Step \ 5 $ for the set of smooth critical points to obtain the last \hl{portions of the} solution curve, i.e. $ B_3 $ and $ B_4 $ trajectories. After doing this, lowest portions of solution curves $ B_1 $ to $ B_{n} $ constitute the solution curve, e.g. curve $ D_iS_1C_2S_3C_4S_5D_f $ in Fig. \ref{fig:Schematic}.

A switching point located on MVC \hl{could either be a smooth} critical point or zero-inertia point. Using the proposed algorithm, depending on the variety of switching points, three general cases may happen which we call them indirect, direct and semi-direct procedures.
When there is no zero-inertia points on MVC, the procedure is indirect. In this case, all of the switching points are to be determined by constructing the portion of MVC which is not covered by initial and final solution curves generated in $ Step \ 1 $.
\hl{The direct} procedure happens when all of the switching points are of the zero-inertia type. In this case, the solution curve is calculated without \hl{the} need to construct any part of MVC.
Semi-direct procedure occurs when switching points are the combination of smooth critical points and zero-inertia points. In this case, the solution is obtained by constructing some portions of MVC which are not covered by \hl{the} generated solution from zero-inertia points in $ Step \ 5 $.
In the next section, to clarify the developed solution algorithm, two examples are given to represent two of these cases.

\section{Numerical Example}
\label{sec:Exm}

In this section, two numerical examples are presented. We consider two different CMMS systems which one of them is planar and the other is three dimensional.
The first example is an example which leads to a semi-direct procedure.
The second example illustrates a case which leads to a completely indirect procedure to solve the time-optimal problem for 3D CMMS that clearly shows the effectiveness and efficiency of the proposed algorithm.

\subsection{Example I: Semi-Direct Procedure}
\label{subsec:Exm_1}

\begin{figure}[!t]
	\centering
	\includegraphics[trim=1cm 3cm 1cm 2cm, clip=true, width=0.9\linewidth]{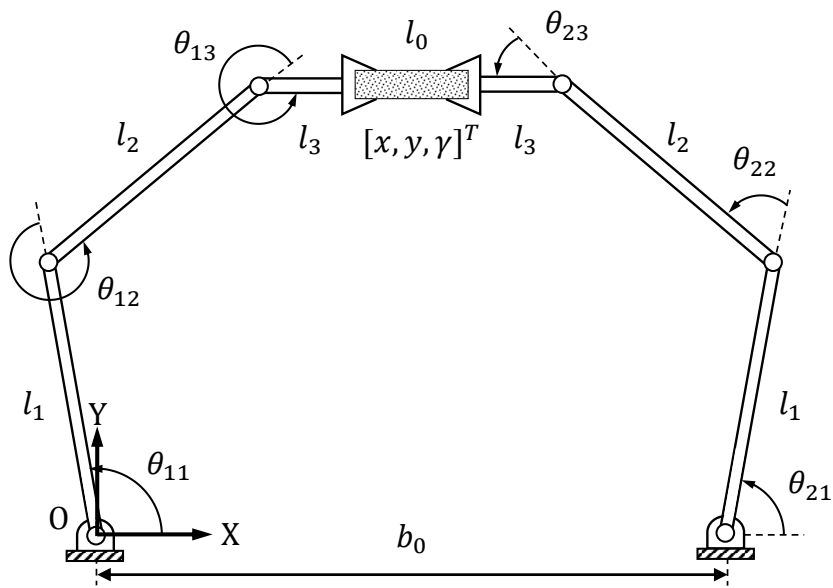}
	\caption{Planar CMMS.}
	\label{fig:planar_CMMS}
\end{figure}

Fig. \ref{fig:planar_CMMS} shows the schematic of a CMMS composed of two planar manipulators handling a payload. Each manipulator has 3 DOFs which rigidly grasped a payload such that no slipping or rotation is possible at contact points. \hl{For this system, the values of $ n $, $ p $ and $ m $ are 6,3 and 6, respectively.} The physical properties of the system are listed in Table \ref{tab:characteristics}. The system is assumed to move the payload on a prescribed path defined as

\begin{table}[!t]
	\renewcommand{\arraystretch}{1.2}
	\caption{Physical characteristics of the System}
	\label{tab:characteristics}
	\centering
	\begin{tabular}{|c|c|c|c|}
		\hline
		\multicolumn{2}{|c|}{length (m)} & \multicolumn{2}{c|}{mass (kg)} \\ \hline
		$ l_0 $                       & 0.2                      & $ m_0 $                      & 1                       \\ \hline
		$ l_1 $                       & 0.5                      & $ m_1 $                      & 1                       \\ \hline
		$ l_2 $                       & 0.6                      & $ m_2 $                      & 1                       \\ \hline
		$ l_3 $                       & 0.3                      & $ m_3 $                      & 0.3                     \\ \hline
		$ b_0 $                       & 1.4                      &                              &                         \\ \hline \hline
		\multicolumn{4}{|c|}{limit torques (N.m)}                                               \\ \hline
		\multicolumn{4}{|c|}{$ \boldsymbol \tau_{limit} = \pm [35, 25, 10, 35, 25, 10]^T $}                                               \\ \hline
	\end{tabular}
\end{table}

\begin{equation}
\label{eq:Path1}
\begin{cases}
x(s)=0.6 (s-0.5)+0.7 \\
y(s)=-3.46 s^5+8.66 s^4-5.77 s^3+0.58 s+0.7 \\
\gamma(s)=-0.45 (s-0.5)
\end{cases}
\end{equation}

\noindent
with $ \dot s_i=4 $ and $ \dot s_f=4 $.

\begin{figure}[!t]
	\centering
	\subfloat[Primary solution curve]{\includegraphics[trim=2cm 6.5cm 2.5cm 0cm, clip=true, width=1\linewidth]{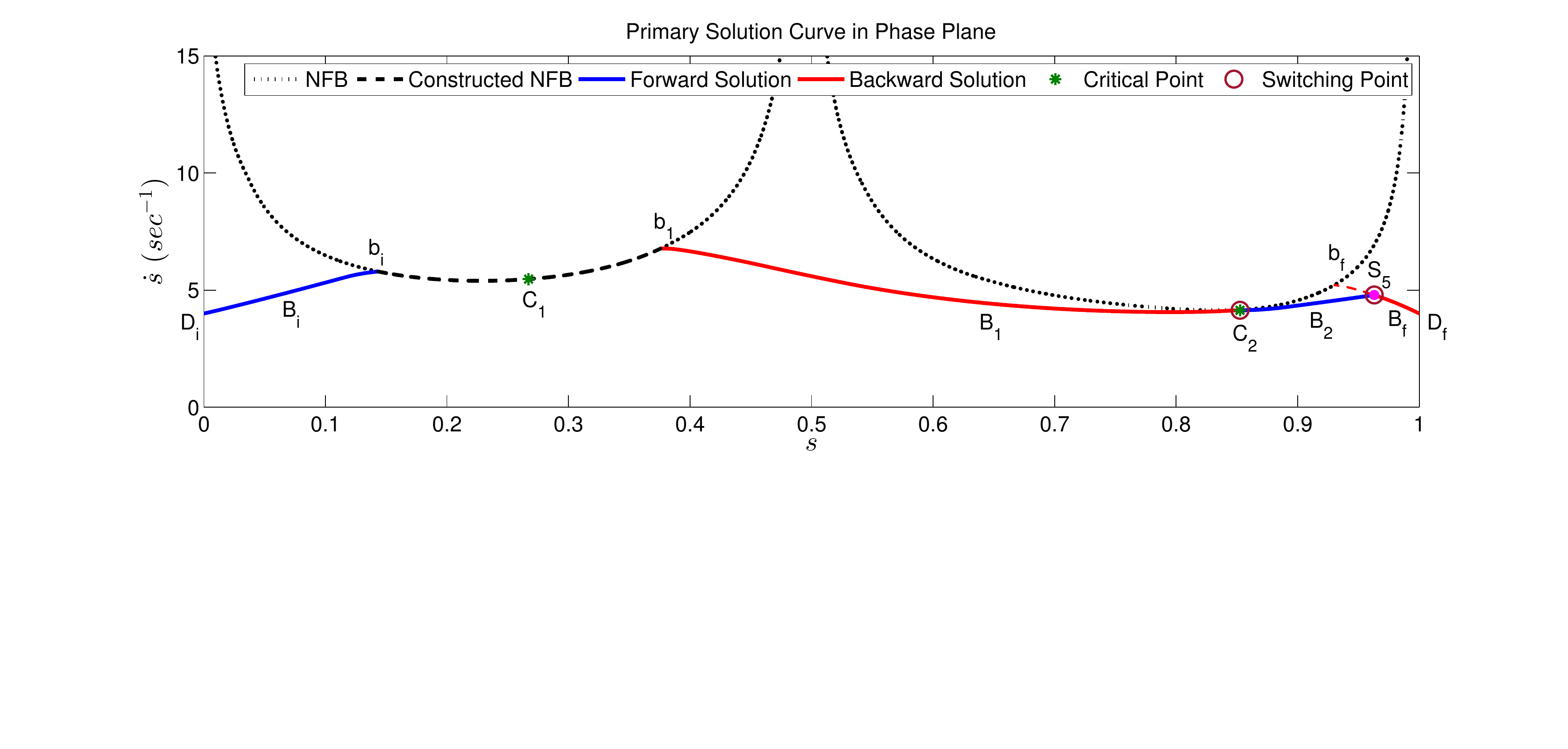}%
		\label{fig:p_solution_curve1}}
	\hfil
	\subfloat[Final solution curve]{\includegraphics[trim=2cm 6.5cm 2.5cm 0cm, clip=true, width=1\linewidth]{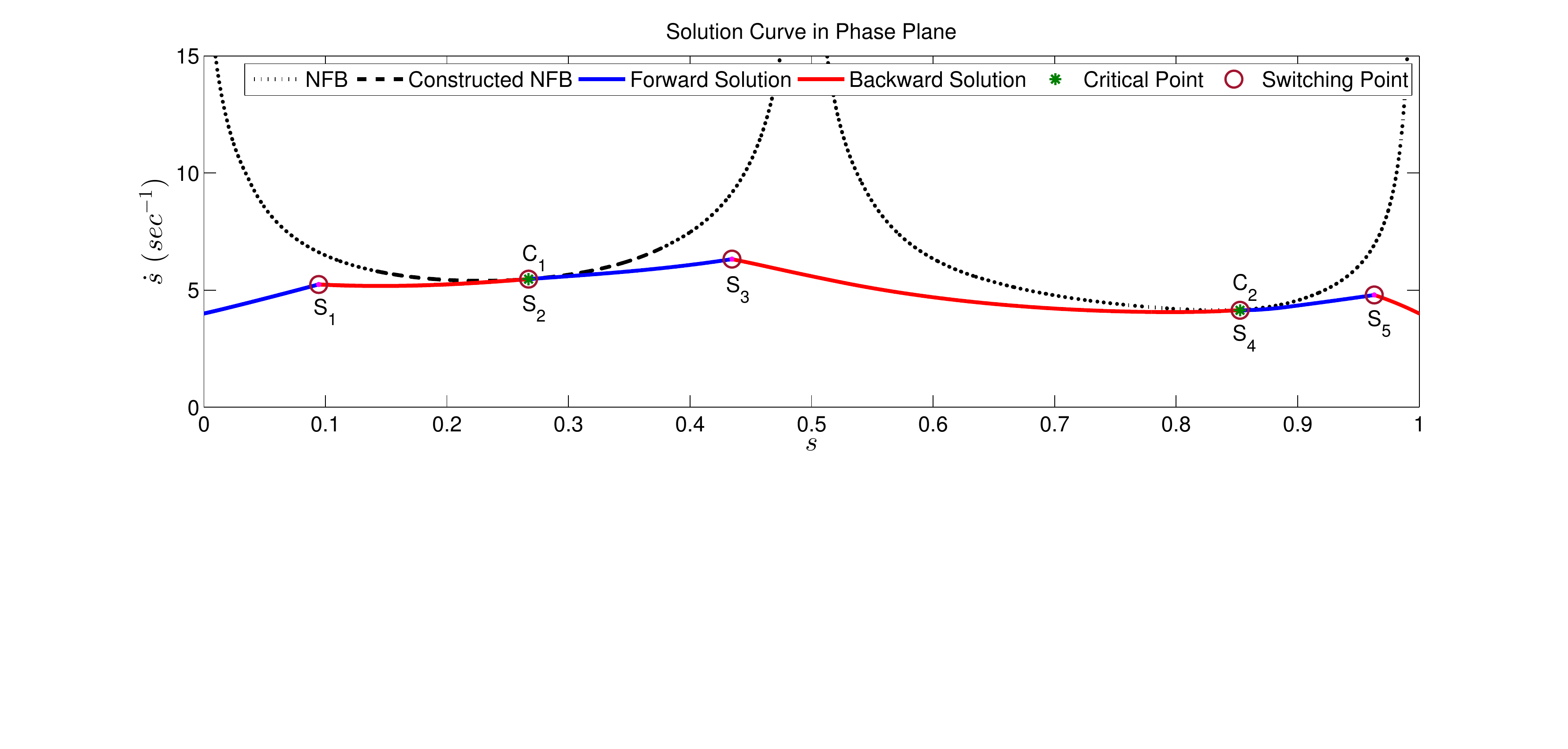}%
		\label{fig:solution_curve1}}
	\caption{Solution curve, critical and switching points calculated for semi-direct procedure (example I).}
\end{figure}

The procedure begins by running SPA strategy in \hl{the} forward direction by starting from \hl{the} initial point, $ D_i $(0, 4), which results branch $ B_i $, as shown in Fig. \ref{fig:p_solution_curve1}. The solution is continued by performing the same strategy in \hl{the} backward direction by starting from the final point, $ D_f $(1, 4), which makes branch $ B_f $. 
Actuator's bounds are violated at $ b_i $(0.1434, 5.7960) and $ b_f $(0.9301, 5.2179) that means the solution enters NFR at these points.
By following the solution algorithm, the only zero-inertia point is calculated as $ C_2 $(0.8526, 4.1395). 
Then, SPA strategy is applied in \hl{the} backward and forward direction by starting from $ C_2 $ that makes branches $ B_1 $ and $ B_2 $. 
When the branch $ B_1 $ reaches to $ b_1 $(0.3758, 6.7744) the actuator's limits are broken. Also, branch $ B_2 $ crosses the curve $ B_f $ at $ S_5 $(0.9630, 4.7920). 

The solution curve is unknown between $ b_i $ and $ b_1 $. Therefore, to complete the solution, a portion of MVC which is not covered by the generated solution curve is constructed. This part of MVC is shown with dashed curve in Fig. \ref{fig:p_solution_curve1}.
\hl{It should be noted} that portions of MVC shown in Fig. \ref{fig:p_solution_curve1} by dotted line are not necessary for calculating of the solution curve and are drawn only for better understanding.

Values of $ \ddot s / \dot s $ and slope of MVC are compared
to categorize points on MVC as sink or source.
The behavior of the MVC changes smoothly from sink to source at $ s_{c_1} = 0.2672 $. 
To obtain the solution curve, we execute SPA strategy by starting from $ C_1 $, forwardly with $ \ddot s_{max} $ and backwardly with $ \ddot s_{min} $ respectively, up to the point where generated solution curve in previous steps is crossed in $ S_1 $(0.0945, 5.2431) and $ S_3 $(0.4345, 6.3232).
The solution curve has been shown in Fig. \ref{fig:solution_curve1}. \hl{Also, the minimum time calculated for moving on this trajectory is $ 0.206 $ sec.}

\begin{figure}[!t]
	\centering
	\includegraphics[trim=4cm 0cm 4cm 0cm, clip=true, width=1\linewidth]{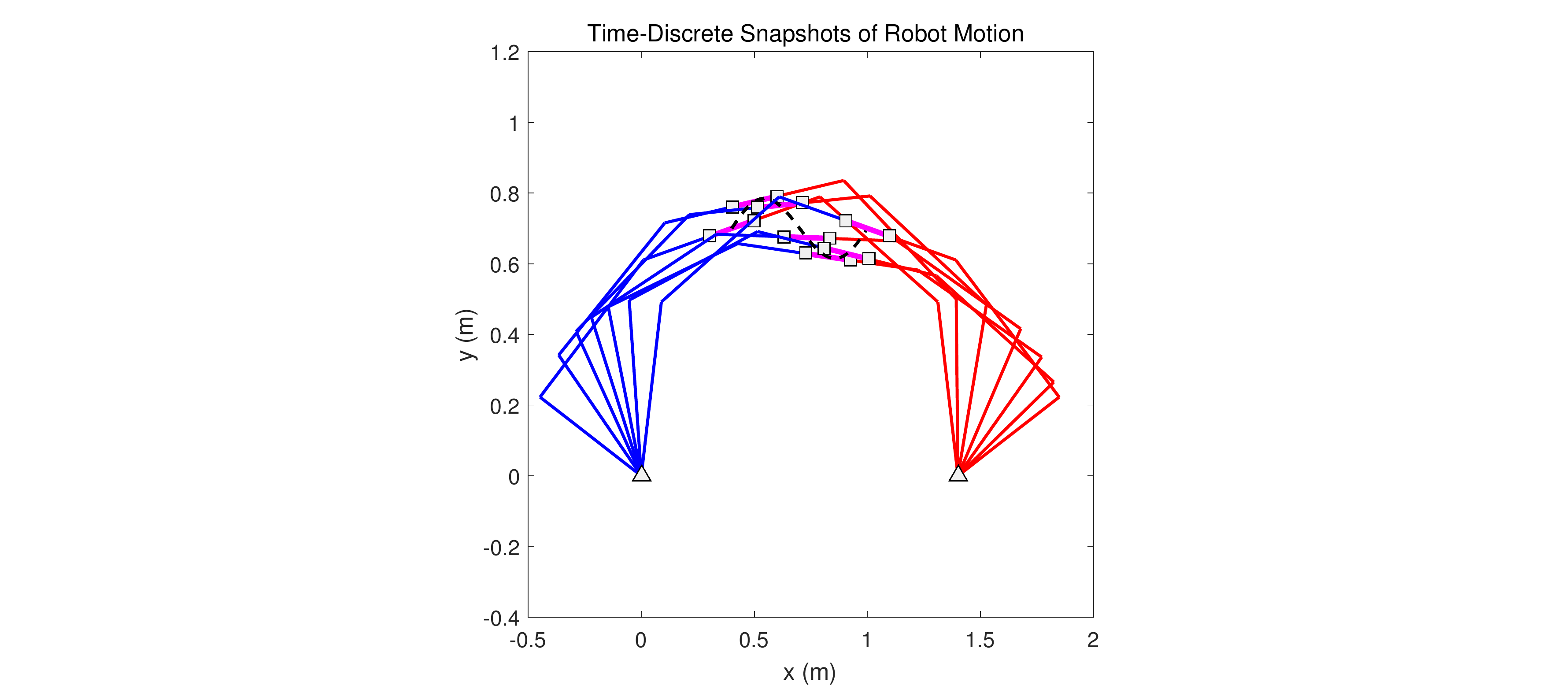}
	\caption{Time-discrete snapshots of robot motion (example I).}
	\label{fig:snapshot}
\end{figure}

\begin{figure}[!t]
	\centering
	\subfloat[Angular position]{\includegraphics[trim=2cm 5.5cm 2.5cm 0cm, clip=true, width=1\linewidth]{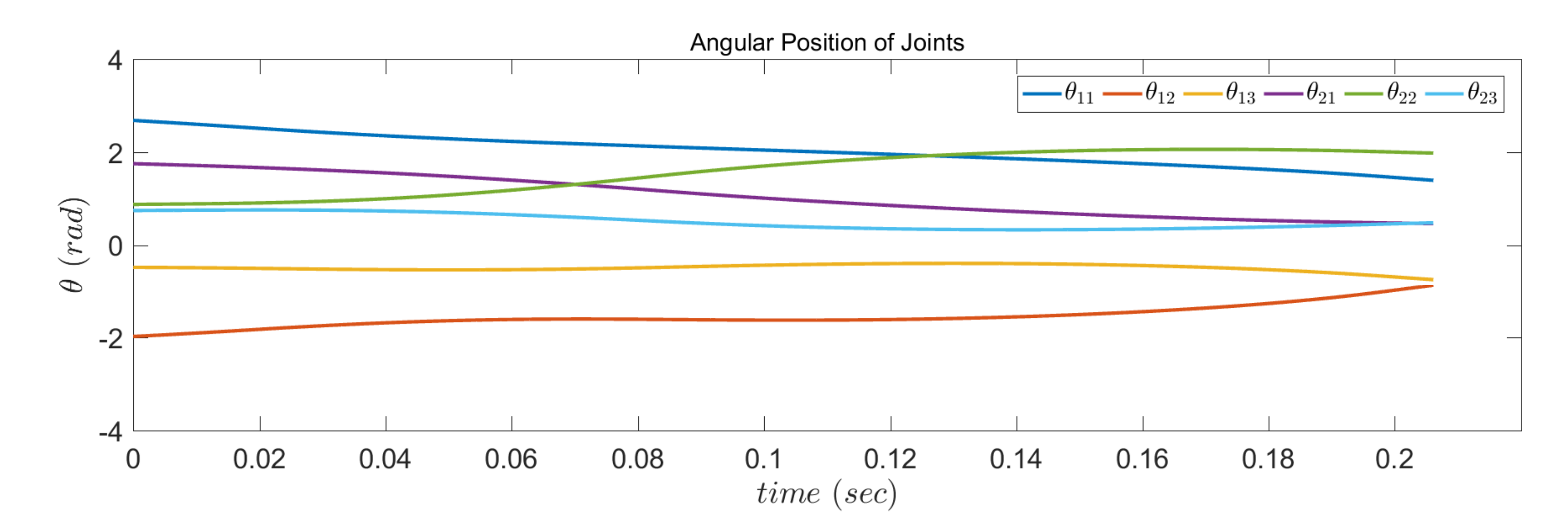}%
		\label{fig:theta}}
	\hfil
	\subfloat[Angular velocity]{\includegraphics[trim=2cm 5.5cm 2.5cm 0cm, clip=true, width=1\linewidth]{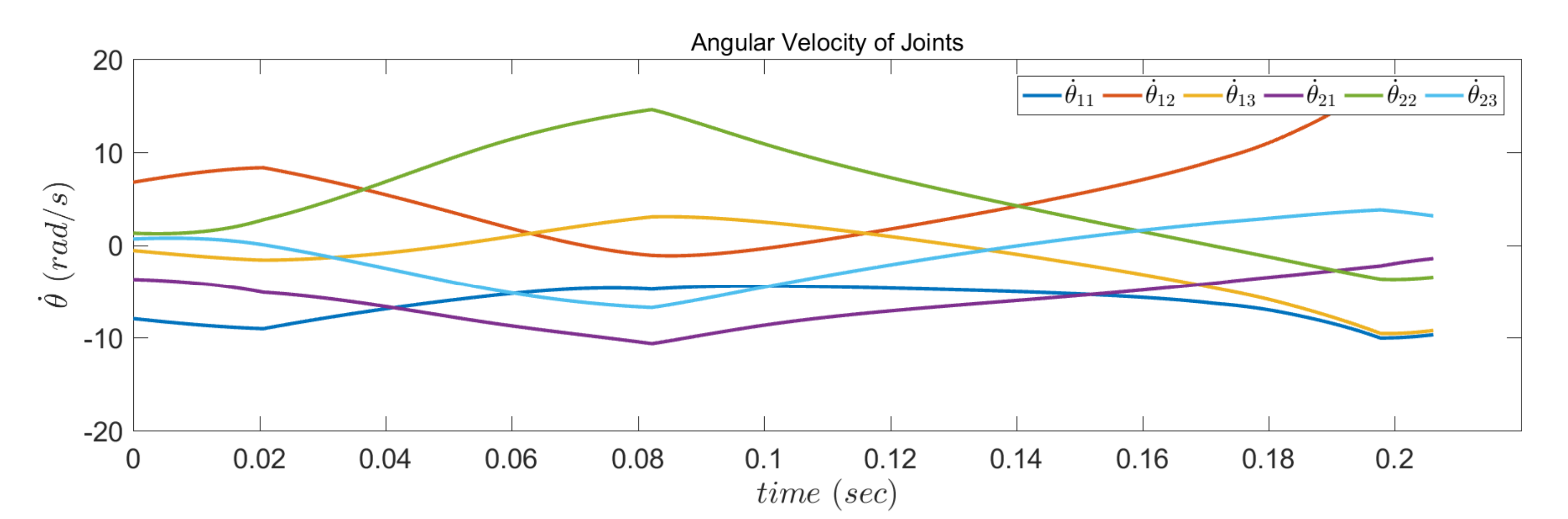}%
		\label{fig:dtheta}}
	\caption{Angular position and velocity of the joints (example I).}
	\label{fig:theta_dtheta}
\end{figure}

\begin{figure}[!t]
	\centering
	\subfloat[Manipulator 1]{\includegraphics[trim=2cm 6.5cm 2.5cm 0cm, clip=true, width=1\linewidth]{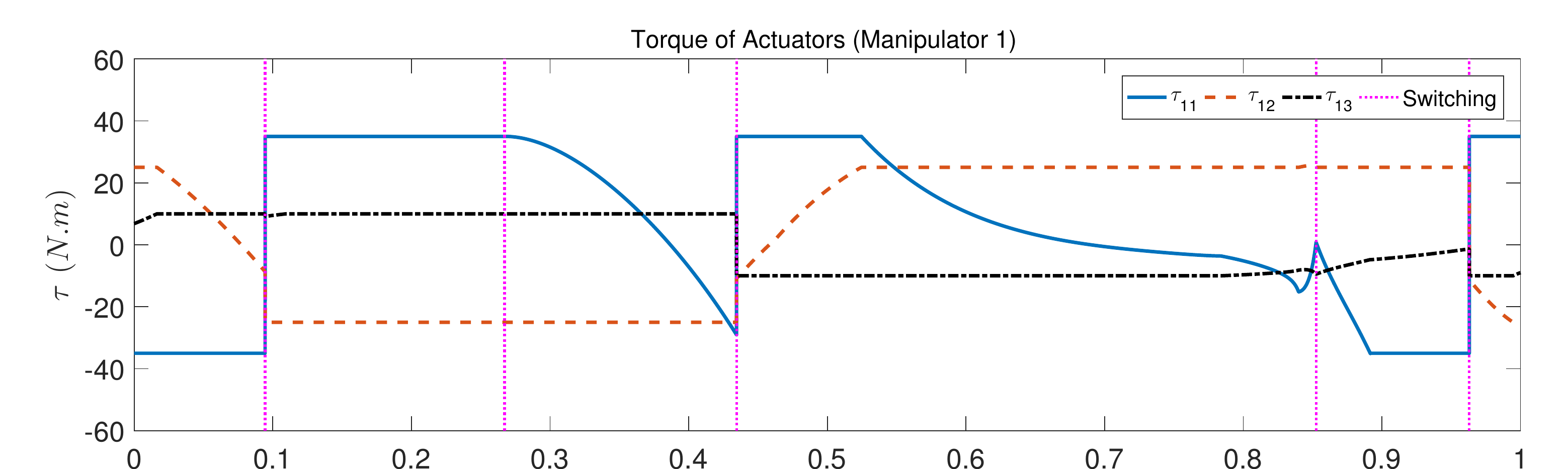}
		\label{fig:tau_manp_1_1}}
	\hfil
	\subfloat[Manipulator 2]{\includegraphics[trim=2cm 6.5cm 2.5cm 0cm, clip=true, width=1\linewidth]{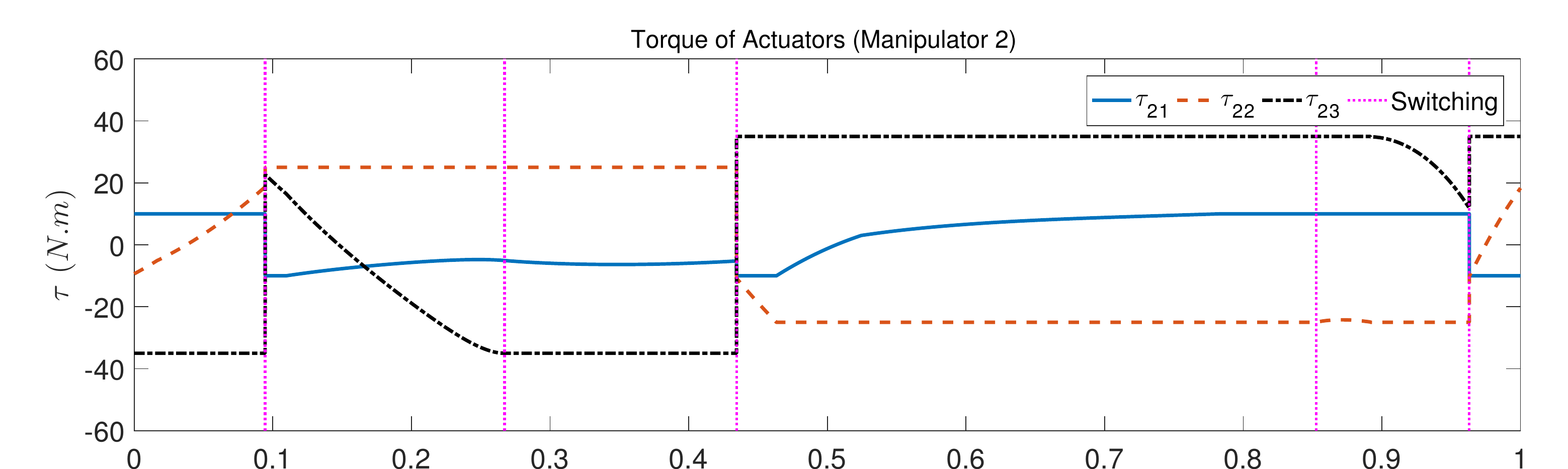}
		\label{fig:tau_manp_2_1}}
	\caption{Actuator torques calculated for semi-direct procedure (example I).}
	\label{fig:tau_manp1}
\end{figure}

The solution curve includes five switching points of $ S_1 $(0.0945, 5.2431), $ S_2 $(0.2672, 5.4703), $ S_3 $(0.4345, 6.3232), $ S_4 $(0.8526, 4.1395) and $ S_5 $(0.9630, 4.7920). Points $ S_1 $, $ S_3 $ and $ S_5 $ are located out of MVC and the acceleration is switched from the maximum value to its minimum. However, two other switching points, $ S_2 $ and $ S_4 $, are placed on MVC and the acceleration is switched from the minimum value to its maximum. 
At both critical points, the acceleration changes continuously from \hl{the} minimum value to its maximum. However, the acceleration is discontinuous at the switching points which are located out of MVC.
Fig. \ref{fig:snapshot} shows time-discrete snapshots of robot motion. The angular position and velocity of the joints are shown in Fig. \ref{fig:theta_dtheta}.
Also, Fig. \ref{fig:tau_manp1} shows the actuator torques of each robot arms while moving on the prescribed path with the minimum time solution. In terms of continuity, the actuator torques has the same behavior of the acceleration at \hl{the} switching points.
Also, the torque limitation has been fulfilled and on each point of the path, at least four actuators are always saturated.

\hl{The solution depends on the system and payload path whether it is obtained by a direct, semi-direct or indirect procedure.  For instance,} it can be easily check that the solution could be obtained directly, i.e. no portion of MVC is needed \hl{to construct,} for payload path defined as

\begin{equation}
\label{eq:Path2}
\begin{cases}
x(s)=0.2 \ cos(2 \pi s)+0.7 \\
y(s)=0.2 \ sin(2 \pi s)+0.62 \\
\gamma(s)=0.7 \ s 
\end{cases}
\end{equation}

\noindent
with $ \dot s_i=0 $ and $ \dot s_f=0 $.

\subsection{Example II: Indirect Procedure}

\begin{figure}[!t]
	\centering
	\includegraphics[width=0.8\linewidth]{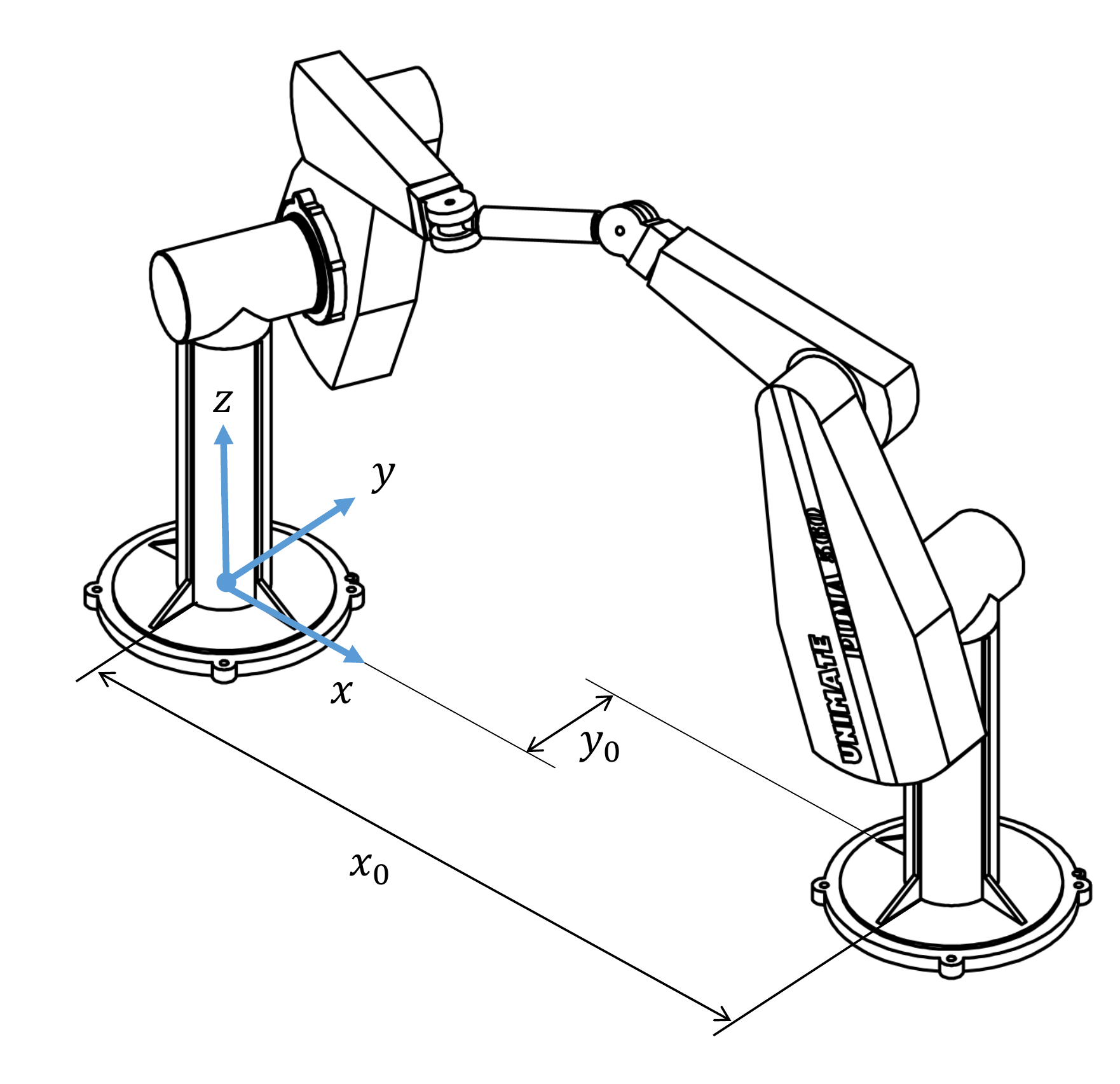}
	\caption{Cooperative multi-manipulator system of two PUMA 560 manipulators.}
	\label{fig:PUMA560}
\end{figure}

The CMMS system considered at this example is composed of two PUMA 560 manipulators which are handling a payload, as shown in Fig. \ref{fig:PUMA560}. 
Each manipulator has 6 DOFs which rigidly grasped a payload such that no slipping or rotation is possible at contact points. \hl{For this system, the values of $ n $, $ p $ and $ m $ are 12, 6 and 12, respectively.} The physical characteristics and Denavit-Hartenberg parameters of PUMA 560 are given in [\citen{Armstrong1986}].
The system is assumed to move the payload on a prescribed path defined as

\begin{equation}
\label{eq:Path3}
\begin{cases}
x(s)=0.2 \ cos(2 \pi s)+0.77 \\
y(s)=0.15 \ sin(2 \pi s)+0.122 \\
z(s)=0.1 \ cos(2 \pi s)-0.1 \\
\alpha(s)=\frac{\pi}{8} \ sin(2 \pi s) \\
\beta(s)=\frac{\pi}{3} \ sin(2 \pi s)+\frac{\pi}{2} \\
\gamma(s)=\frac{\pi}{12} \ sin(2 \pi s)+\frac{\pi}{6} \\
\end{cases}
\end{equation}

\noindent
with $ \dot s_i=1 $ and $ \dot s_f=1 $.

\begin{figure}[!t]
	\centering
	\subfloat[Primary solution curve]{\includegraphics[trim=2cm 6.5cm 2.5cm 0cm, clip=true, width=1\linewidth]{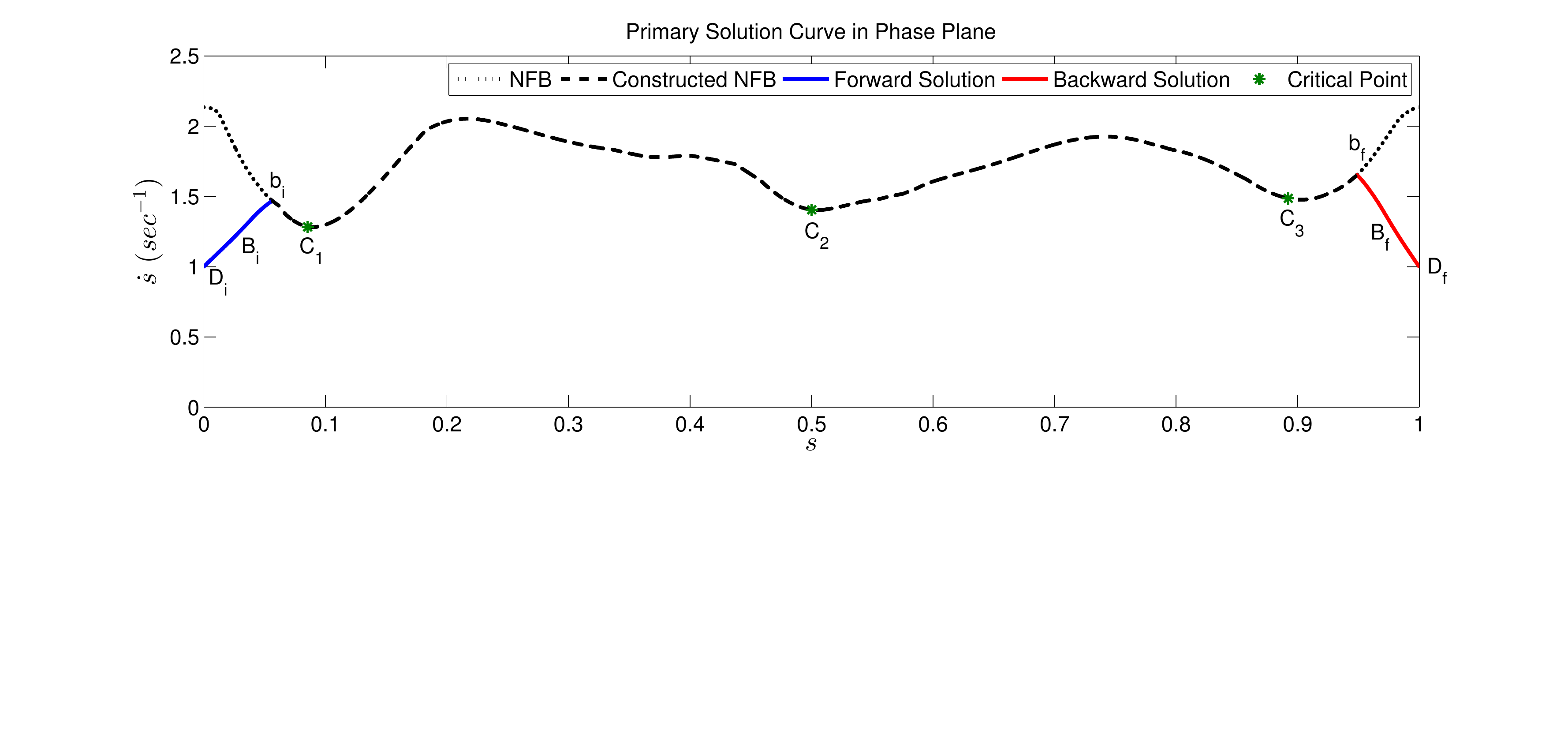}%
		\label{fig:p_solution_curve3}}
	\hfil
	\subfloat[Final solution curve]{\includegraphics[trim=2cm 6.5cm 2.5cm 0cm, clip=true, width=1\linewidth]{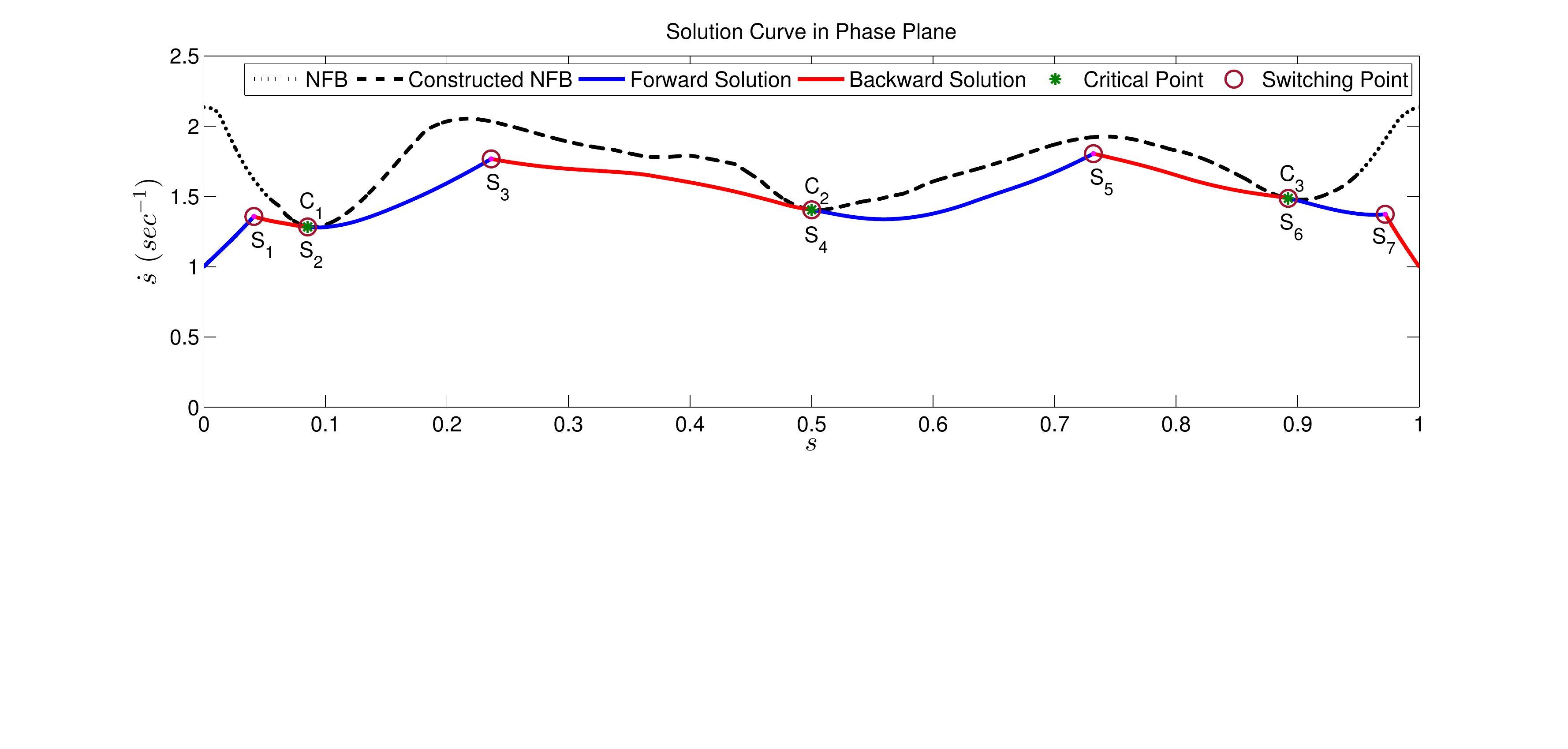}%
		\label{fig:solution_curve3}}
	\caption{Solution curve, critical and switching points calculated for indirect procedure (example II).}
\end{figure}

The procedure begins by running SPA strategy in \hl{the} forward and backward direction from $ D_i $(0, 1) and $ D_f $(1, 1) respectively which results in branches $ B_i $ and $ B_f $, as shown in Fig. \ref{fig:p_solution_curve3}. \hl{Actuator's bounds have been violated} at $ b_i $(0.0564, 1.4694) and $ b_f $(0.9488, 1.6534) that means the solution enters to NFR at these points.
In the range of $ b_ib_f $, there is no zero-inertia point of the sink-source type. Thereby to generate the solution, a portion of MVC placed between $ b_i $ and $ b_f $ is constructed. This part of MVC is shown with dashed curve in Fig. \ref{fig:p_solution_curve3}.

\begin{figure}[!t]
	\centering
	\subfloat[Manipulator 1]{\includegraphics[trim=2cm 5.5cm 2.5cm 0cm, clip=true, width=1\linewidth]{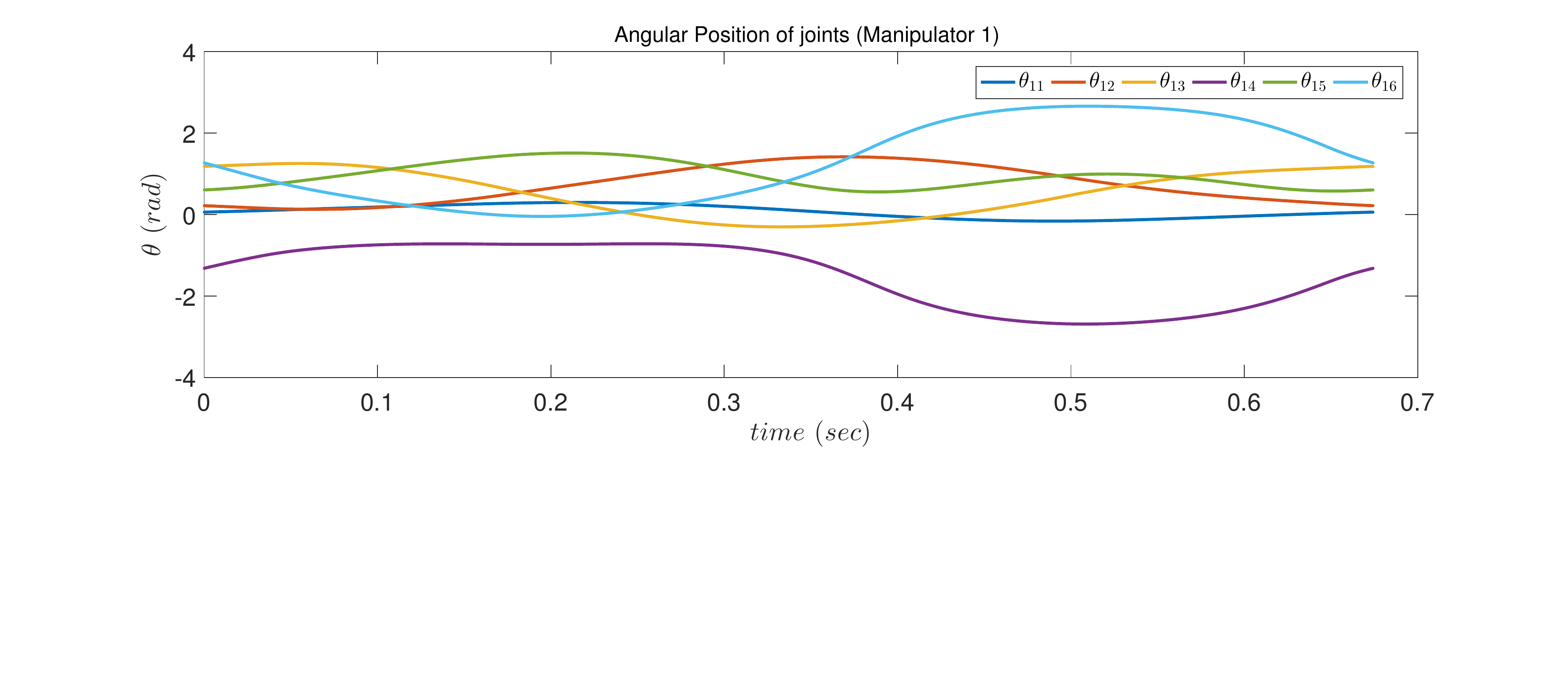}%
		\label{fig:theta1}}
	\hfil
	\subfloat[Manipulator 2]{\includegraphics[trim=2cm 5.5cm 2.5cm 0cm, clip=true, width=1\linewidth]{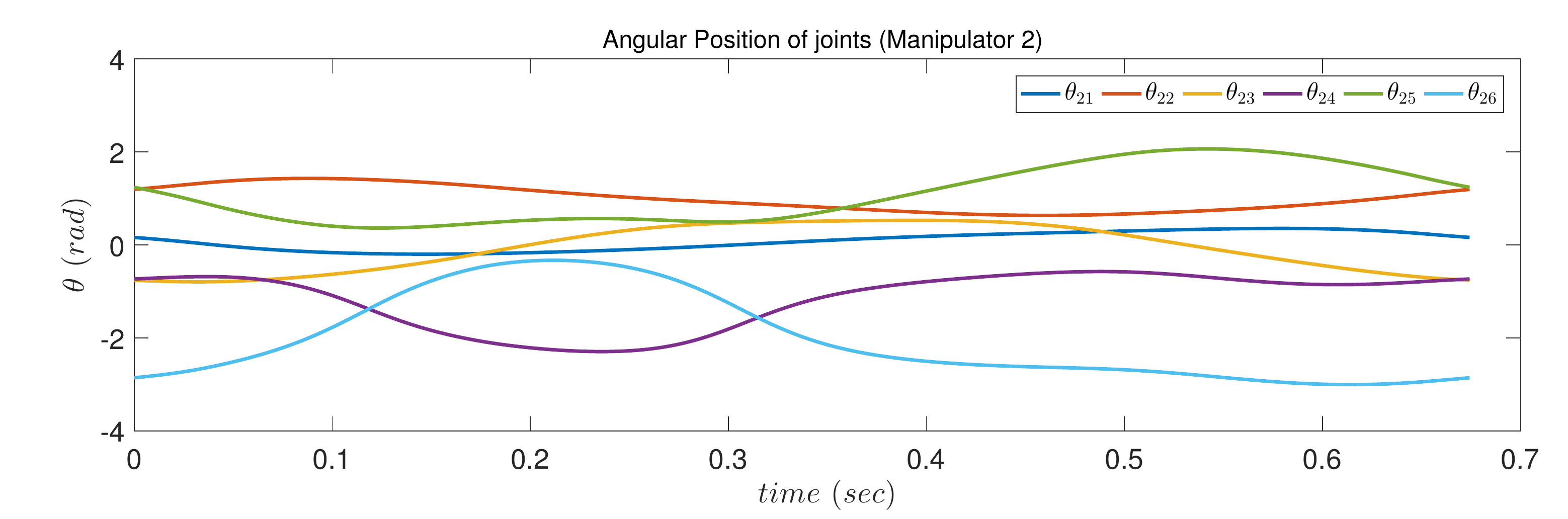}%
		\label{fig:theta2}}
	\caption{Angular position of the joints (example II).}
	\label{fig:theta12}
\end{figure}

\begin{figure}[!t]
	\centering
	\subfloat[Manipulator 1]{\includegraphics[trim=2cm 5.5cm 2.5cm 0cm, clip=true, width=1\linewidth]{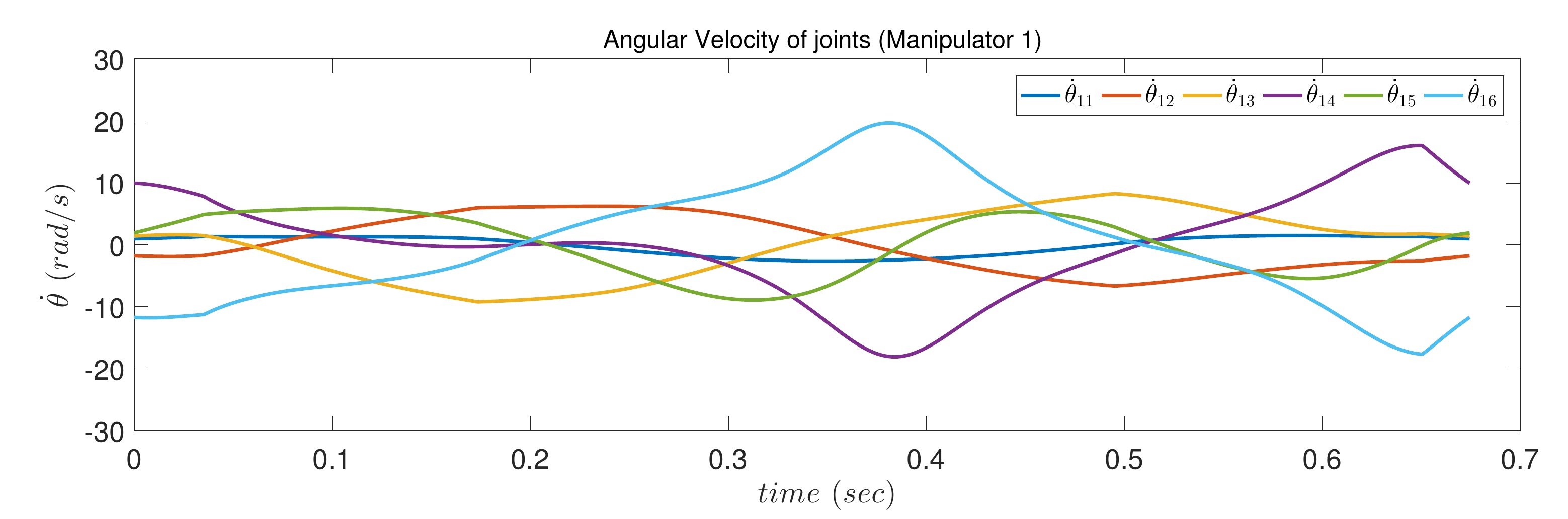}%
		\label{fig:dtheta1}}
	\hfil
	\subfloat[Manipulator 2]{\includegraphics[trim=2cm 5.5cm 2.5cm 0cm, clip=true, width=1\linewidth]{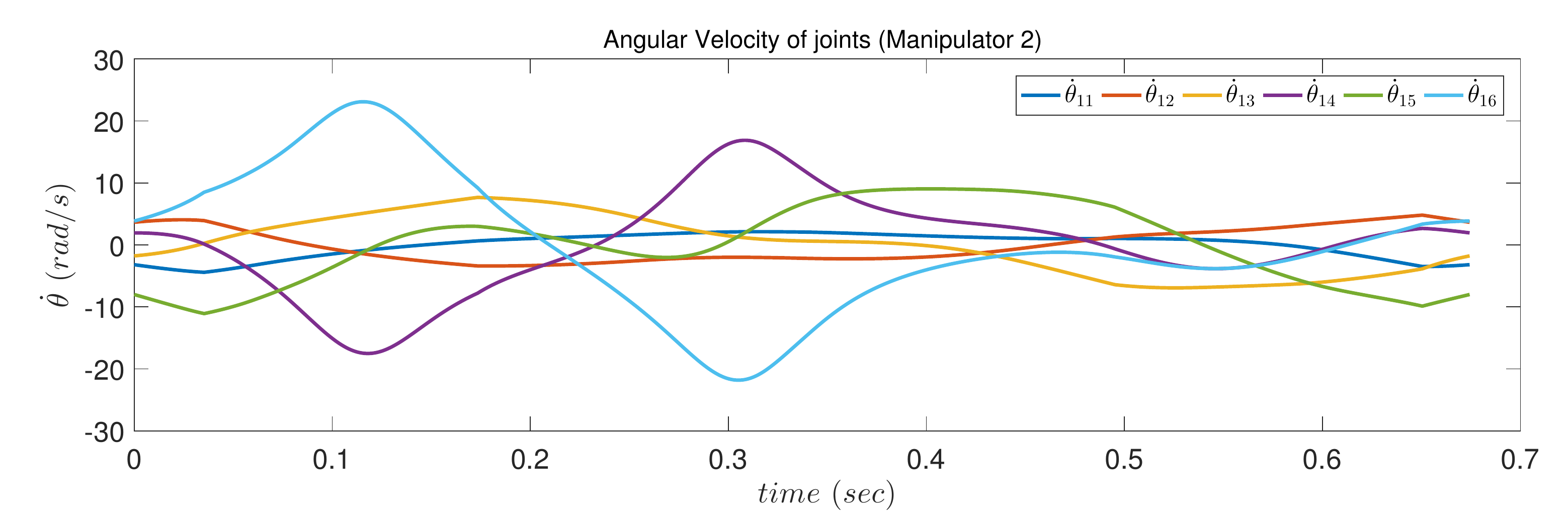}%
		\label{fig:dtheta2}}
	\caption{Angular velocity of the joints (example II).}
	\label{fig:dtheta12}
\end{figure}

\begin{figure}[!t]
	\centering
	\subfloat[Manipulator 1]{\includegraphics[trim=2cm 6.5cm 2.5cm 0cm, clip=true, width=1\linewidth]{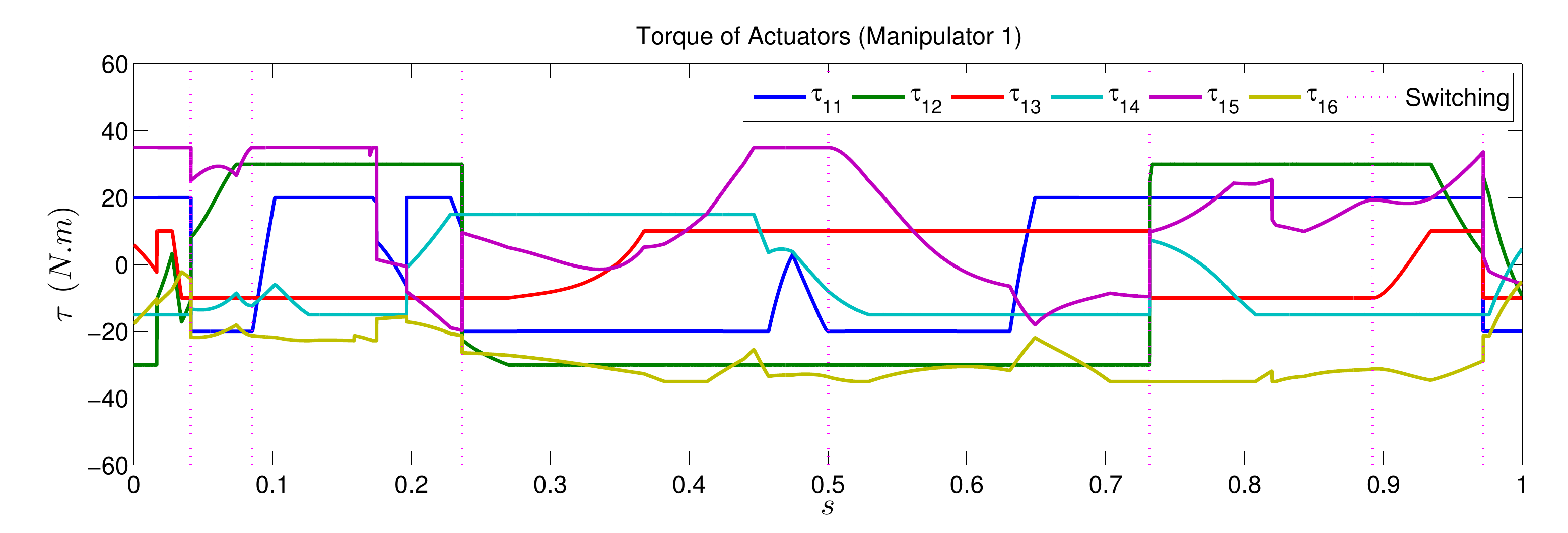}%
		\label{fig:tau_manp_1_3}}
	\hfil
	\subfloat[Manipulator 2]{\includegraphics[trim=2cm 6.5cm 2.5cm 0cm, clip=true, width=1\linewidth]{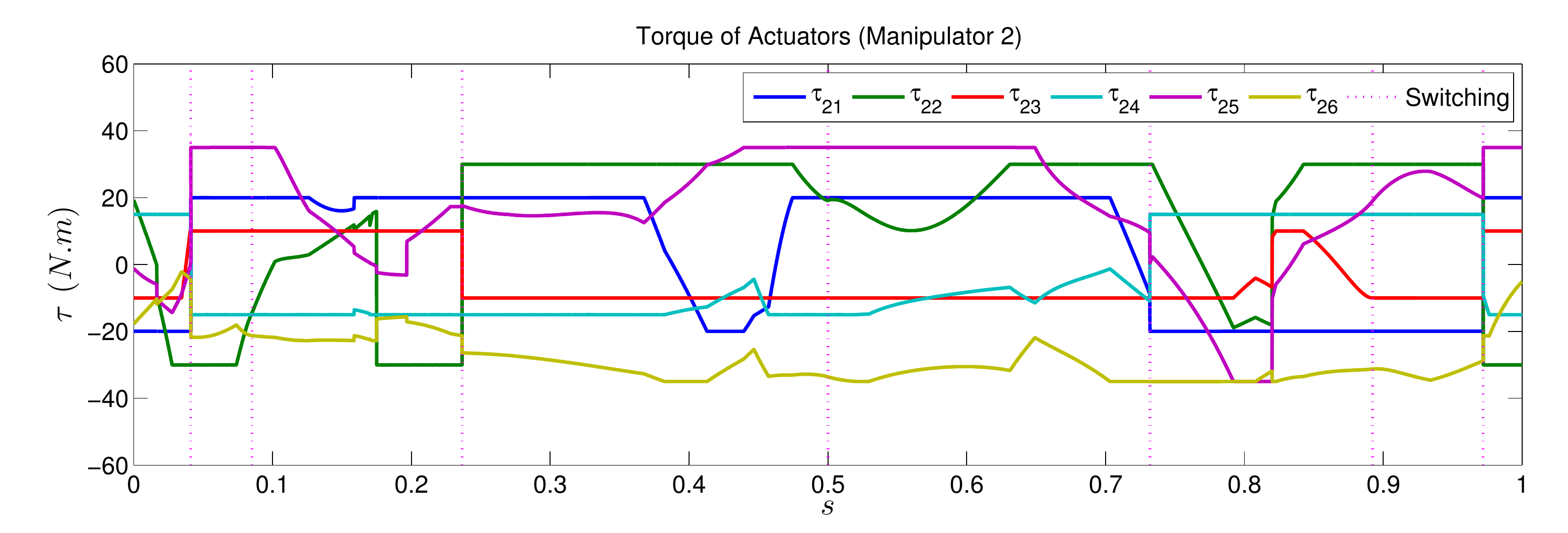}%
		\label{fig:tau_manp_2_3}}
	\caption{Actuator torques calculated for indirect procedure (example II).}
	\label{fig:tau_manp3}
\end{figure}

The values of $ \ddot s / \dot s $ and slope of MVC are compared to determine \hl{the} smooth critical points.  
The behavior of the MVC changes smoothly from sink to source at points $ s_{c_1} = 0.0854 $, $ s_{c_2} = 0.5001 $ and $ s_{c_3} = 0.8922 $.
The solution curves are established by performing SPA strategy in \hl{the} backward and forward direction by starting from points $ C_1 $, $ C_2 $ and $ C_3 $ respectively. The generated curves cover whole domain of $ s $ from 0 to 1 and the solution algorithm ends. 
The solution curve consists of seven switching points $ S_1 $(0.0412, 1.3573), $ S_2 $(0.0854, 1.2833), $ S_3 $(0.2365, 1.7672), $ S_4 $(0.5001, 1.4047), $ S_5 $(0.7319, 1.8045), $ S_6 $(0.8922, 1.4868) and $ S_7 $(0.9720, 1.3734) as shown in Fig. \ref{fig:solution_curve3}. The minimum time calculated, also, amounts to $ 0.674 $ sec. The angular position and velocity of the joints are shown in Fig. \ref{fig:theta12} and Fig. \ref{fig:dtheta12}, respectively.
The actuator torques of each robot arms for the minimum time solution is shown in Fig. \ref{fig:tau_manp3}. The torque limitation has been fulfilled on each point of the path and at least seven actuators are always located in their saturation bound.

\hl{The importance of using an efficient algorithm to solve the minimum time problem is more clear in $ Example \ II $.} It can be investigated in two different aspects. First, the CMMS system which is composed of two PUMA 560 manipulators is more complicated than planar robot considered in  the previous example.
Second, for the path prescribed in (\ref{eq:Path3}), there is no zero-inertia point, therefore, to determine the critical points, a large portion of MVC should be constructed.

\section{Efficiency Evaluation}
\label{sec:EE}

\begin{table*}
	\renewcommand{\arraystretch}{1.2}
	\centering
	\caption{Evaluating the efficiency of proposed time-optimal algorithm by comparison with reference [\citen{Pham2015}]}
	\label{tab:eff}
	\begin{tabular}{|c|c|c|c|c|}
		\hline
		index                                                                        & \multicolumn{4}{c|}{time consumption (sec)}           \\ \hline
		example                                                                      & \multicolumn{2}{c|}{I}    & \multicolumn{2}{c|}{II}   \\ \hline
		method                                                                       & ref [\citen{Pham2015}]        & proposed     & ref [\citen{Pham2015}]      & proposed      \\ \hline \hline
		MVC                                                                          & 5.4        & 1.2          & 696       & 594           \\ \hline
		solution curve                                                               & 12.6       & 5.6          & 1608      & 402           \\ \hline
		time-optimal path                                                            & 18         & 6.8          & 2304      & 996           \\ \hline
		\begin{tabular}[c]{@{}c@{}}reduction of\\ time consumption (\%)\end{tabular} & \multicolumn{2}{c|}{62.2} & \multicolumn{2}{c|}{56.8} \\ \hline
	\end{tabular}
\end{table*}

In this section, the efficiency of the proposed algorithm is compared \hl{with one} of the most recent method. We implemented the presented method in [\citen{Pham2015}] for the numerical examples solved in section \ref{sec:Exm}. Both methods are run on Matlab/Simulink 9.0 enabling parallel computing on a PC (Intel i7-5930 CPU and 24 GB Memory). The forward and backward integrations performed with 1 $ ms $ as step time. The MVC is also constructed \hl{by considering 1000 points along the path.} The computation time for calculating two major parts of the time-optimal path, i.e. MVC and solution curve, are given in Table \ref{tab:eff} for two methods as well as total computation time. The results show that the computation time are reduced 62.2\% and 56.8\% for examples I and II, respectively. The time consumed to calculate the MVC for the examples is reduced which is due to skip of the whole construction of MVC by implementing our algorithm. 
Using the SPA algorithm, the computation time for calculating solution curve is substantially reduced by looking for saturated actuators in only a limited number of points, instead of determining the saturated actuators in each point which is performed in [\citen{Pham2015}].
For example II, the equations of motion become more complicated which results a large computation cost to find \hl{the} time-optimal path. By implementing the proposed method, the
computation time for applying the method presented in [\citen{Pham2015}], which \hl{was} 2304 $ sec $,  is drastically reduced to 996 $ sec $. This example clearly shows the efficiency of our algorithm \hl{by saving of} 1308 $ sec $.

We also implemented our method for the model of biped robot reported  by  Liu \textit{et al.} \cite{Liu2013}. They used parametric optimization method to find periodic steady state trajectory  for a biped walking at \hl{a specified} speed of 1.8 $ km/h $, and tried to re-optimize the resulted trajectory \hl{by} using dynamic differential programming method. As proposed in our previous study \cite{Sadigh2013},  we formulated the minimum time problem with stability and non-slip conditions along with actuator's limits expressed as some inequality constraints. \hl{In addition, certain kinematic constraints are considered in terms of hip joint position that ensure an acceptable walking pattern.} By solving this problem with the method presented \hl{in this paper, the maximum achievable speed is 5.6 $ km/h $.}

\section{Conclusion}
The problem of finding the time-optimal trajectory of a cooperative multi-manipulator system with redundant actuators moving on a specified path is studied in this paper. \hl{The focus of this study is to reduce the} computational cost of the solutions obtained based on phase plane analysis. Two important issues which are sources of large computational costs were considered and resolved. The first \hl{problem} is how to compute minimum or maximum acceleration along the path. Using the presented method, instead of determining the saturated actuators in each point which is computationally expensive, we can look for saturated actuators \hl{only in} a limited number of points. The second problem is \hl{the} calculation of the switching points which is the most difficult part of the time-optimal solution. 
\hl{This paper provides a comprehensive study on this subject and examines all possible situations.} 
Several examples are \hl{provided to show} different possibilities which may be encountered, \hl{and how} they are solved by means of the proposed method.

\end{document}